\documentclass[11pt]{amsart}

\usepackage{hyperref}

\usepackage{amssymb, amsmath, amsthm,amsrefs}

\usepackage{latexsym}
\usepackage{color}
\usepackage{exscale}
\usepackage{fullpage}
\usepackage{bm, bbm, mathrsfs}

\newtheorem{theorem}{Theorem}[section]
\newtheorem{lemma}[theorem]{Lemma}
\newtheorem{proposition}[theorem]{Proposition}

\newtheorem{problem}{Problem}
\newtheorem{definition}{Definition}

\newtheorem{example}[theorem]{Example}

\newtheorem{remark}{Remark}[section]

\def\s0{{ s_0}}

\def\ts0{{\tilde s_0}}
\def\TA{{\tilde A}}
\def\tor{{\bf T}}
\def\cur{{\bf R}}

\def\TS{\textstyle}

\def\eq#1{(\ref{#1})}
\def\nn{\nonumber}

\def\({\left(\begin{array}{cccccc}}
\def\){\end{array}\right)}

\def\bes{\begin{eqnarray}}
\def\ees{\end{eqnarray}}

\def\dbyd#1{\TS{\frac{\partial}{\partial#1}}}

\def\point{{\mathbf p}}

\newcommand{\del}{\partial}
\newcommand{\pd}[2]{\frac{\partial{#1}}{\partial{#2} }}
\newcommand{\beq}{\begin{equation}}
\newcommand{\eeq}{\end{equation}}
\newcommand{\bea}{\begin{eqnarray}}
\newcommand{\eea}{\end{eqnarray}}
\newcommand{\beann}{\begin{eqnarray*}}
\newcommand{\eeann}{\end{eqnarray*}}

\newcommand{\lam}{\ensuremath{\lambda}}

\newcommand{\RR}{\mathbb{R}}

\newcommand{\ti}{\tilde}
\newcommand{\ka}{\ensuremath{\kappa}}
\newcommand{\ha}{\ensuremath{h}}
\newcommand{\tha}{\ensuremath{ {\tilde h}}}

\newcommand\br[1]{\overline{{#1}}}
\newcommand{\pf}{\begin{proof}}
\newcommand{\foorp}{\end{proof}}

\DeclareMathOperator{\spa}{span}

\DeclareMathOperator{\rank}{rank}

\DeclareMathOperator{\codim}{codim}
\DeclareMathOperator{\diag}{diag}

\numberwithin{equation}{section}
\begin{document}
\title{Systems of hyperbolic conservation laws with prescribed eigencurves}
\author{Helge Kristian Jenssen }\address{ H.~K.~Jenssen, Department of
Mathematics, Penn State University,
University Park, State College, PA 16802, USA ({\tt
jenssen@math.psu.edu}).}\thanks{ H.~K.~Jenssen was partially supported by NSF grant
DMS-0539549 and worked  on this paper during his participation in the  international research program on
Nonlinear Partial Differential Equations at the Centre for Advanced
Study at the Norwegian Academy of Science and Letters in Oslo during
the academic year 2008Ð-09.}
\author{Irina A. Kogan}
\address{Irina A. Kogan, Department of Mathematics,
North Carolina State University, Raleigh, NC, 27695, USA  ({\tt
iakogan@ncsu.edu}).}\thanks{ I.~A.~Kogan was partially supported by NSF grant
CCF-0728801}

\date{\today}

\begin{abstract}
	We study the problem of constructing systems of hyperbolic conservation laws 
	in one space dimension with prescribed eigencurves, i.e.\ the eigenvector 
	fields of the Jacobian of the flux are given. We formulate this as a typically overdetermined
	system of equations for the eigenvalues-to-be. Equivalent formulations in terms of differential 
	and algebraic-differential equations are considered. The resulting equations  are then
	analyzed using appropriate integrability theorems (Frobenius, Darboux and Cartan-K\" ahler). 
	We give a complete analysis of the possible scenarios, including examples, for systems of three 
	equations. As an application we characterize conservative systems with the same eigencurves
	as the Euler system for 1-dimensional compressible gas dynamics. 
	The case of general rich systems of any size (i.e.\ when the given eigenvector fields are 
	pairwise in involution; this includes all systems of two equations) is completely resolved and we consider 
	various examples in this class. 	
\end{abstract}
\maketitle



%
%
\section{Introduction}\label{intro}
\noindent
Consider a system of $n$ conservation laws  
in one space dimension written in canonical form, 
\beq\label{claw}
	u_t+f(u)_x = 0\,.
\eeq
Here the unknown state $u=u(t,x)\in\mathbb R^n$ is assumed to vary 
over some open subset $\Omega\subset\mathbb R^n$ and the flux $f=(f^1,\dots,f^n)^T$ is a 
nonlinear map from $\Omega$ into $\mathbb R^n$. 
The eigenvalues and eigenvectors of the Jacobian 
matrix $Df(u)$ provide information that is used to solve the Cauchy
problem for \eq{claw}. In particular, the geometric properties 
of the integral curves of eigenvector fields of $Df$ play a key role. 
Together with the so-called Hugoniot curves (see below for definition) these
form {\em wave curves} in $u$-space that are used to build solutions to \eq{claw}.

There is a well developed theory for the Cauchy problem for a large 
class of  systems \eq{claw} in the near equilibrium regime 
$u(0,x)\approx$ constant. Glimm \cite{gl} established 
global-in-time existence of weak solutions for data with sufficiently 
small total variation. By now there are several methods available to obtain 
these solutions: Glimm's original random choice method  \cite{gl},
Liu's deterministic version \cite{liu1}, the wave-front tracking schemes of 
DiPerna \cite{dip1}, Bressan \cite{br1} and Risebro \cite{ri}, and the 
wave tracing methods of Bianchini and Bressan \cite{bb}, \cite{bi1}, \cite{bi2}.

Single equations and systems of two equations enjoy particular properties that yield global 
existence beyond the perturbative regime, \cite{kru}, \cite{gla}, \cite{ni}, \cite{nism}, 
\cite{tar}, \cite{dip2}, \cite{chen}. On the other hand there is no general 
existence result available for ``large" data when the system has three or 
more equations (e.g.\ the Euler system for compressible gas dynamics). 
Indeed, examples of blowup in total variation and/or sup norm are known 
even for genuinely nonlinear and  strictly hyperbolic systems, \cite{jmr}, 
\cite{yo}, \cite{je}, \cite{jy}, \cite{bj}. These examples show that the global 
behavior of the wave curves in $u$-space is related to finite time blow up, 
\cite{yo}, \cite{je}, \cite{jy}, \cite{bj}. It is of interest to have more examples 
of this type. In particular one would like to know if a physical system can 
display similar behavior\footnote{This refers to the strictly hyperbolic regime. 
It is well-known that failure of strict hyperbolicity can cause singular behavior 
(blowup in total variation or $L^\infty$), even in physical systems.
See \cite{te}, \cite{sever} and references therein.}. In searching for systems 
whose wave curves have special properties one is naturally lead to ask 
what freedom one has in {\em prescribing} such curves.

In this paper we consider the situation where we are given a frame of $n$ linearly 
independent vector fields and their integral curves. 
We then ask if there are any systems of the form \eq{claw} with the property that 
the given vector fields are the eigenvector fields of $Df(u)$. 

It turns out that there is always a one-parameter family of trivial solutions, and that
there are cases where no non-trivial solutions exist. More generally we are interested
in knowing ``how many" solutions there are, in terms of how many constants and functions
appear in a general solution. The basic tool for this sort of questions is provided by the 
integrability theorems for overdetermined systems of PDEs  such as Frobenius, Darboux and 
Cartan-K\"ahler theorems. The latter requires re-writing the system of PDE's as an exterior 
differential systems (EDSs) \cite{bcggg}, \cite{il}.

The case of $2\times 2$-systems was treated in \cite{daf1}; it is also covered by 
the analysis of rich systems in Section \ref{rich}.
On the other hand we shall see that already for systems with three equations there
are several possible scenarios. 
Before giving a precise formulation of the problem we review some relevant background material.

{\bf Notation:}
	We denote the $(i,j)$-entry (i.e., the element in the $i$th row and the $j$th column) 
	of an $m\times n$-matrix $A$ by $A^i_j$. Superscript $\,{}^T$ denotes transpose.
	Summation convention is {\em not} used.

\subsection{Conservation laws in one space dimension}\label{prelim-cl}
We consider {\em hyperbolic} systems of conservation laws in one spatial dimension \eq{claw}, 
i.e.\ the Jacobian matrix $Df(u)$ is diagonalizable over $\mathbb R$ at each state 
$u \in\Omega$. The system is {\em strictly} hyperbolic in $\Omega$ provided the eigenvalues 
$\lam^i(u)$ of $Df(u)$ are real and distinct:
\beq\label{strict_hyp}
	\lam^1(u)<\cdots<\lam^n(u)\,,\qquad\qquad\forall\,  u\in\Omega\,.
\eeq
Let us for now fix a choice of the associated right and left eigenvectors 
$R_i(u)$ and $L^i(u)$ of $Df(u)$. These are considered as column and row 
vectors of functions, respectively, and we write
\[R_i(u)=\big[R^1_i(u),\dots,R^n_i(u)\big]^T\,,\qquad\qquad 
L^i(u)=\big[L^i_1(u),\dots,L^i_n(u)\big]\,.\] 
We refer to the $R_i(u)$ as the {\em eigenfields} and their integral curves in $u$-space
as {\em eigencurves}. Diagonalizing $Df$ we have
\beq\label{diag}
	Df(u)=R(u)\, \Lambda(u)\, L(u)\,,
\eeq
where 
\[R(u)=[R_1(u)\,|\,\cdots\,|\, R_n(u)]\,,\quad 
\Lambda(u)=\mbox{diag}[\lam^1(u)\dots\lam^n(u)]\,, \]
and
\[L(u)=R(u)^{-1}=\left[\begin{array}{c} L^1(u)\\ \hline \vdots\\\hline L^n(u)\end{array}\right]\,.\]
Note that in setting $L(u)=R(u)^{-1}$ we have introduced the normalization 
$R_i(u)\cdot L^j(u)=\delta_{i}^{j}$ (Kronecker delta).

Next consider the initial value problem for \eq{claw} where the data consist of two 
constant states separated by a jump at $x=0$,
\beq\label{ri}
	u_0(x)=\left\{\begin{array}{ll}
	u_-\,,& x<0\\
	u_+\,,& x>0\,.
	\end{array}\right.
\eeq
This is the so-called Riemann problem and its solution serves as a building block for
more general solutions, \cite{lax}, \cite{gl}. For a sufficiently smooth flux $f$ whose characteristic fields are 
genuinely nonlinear or linearly degenerate in the sense of Lax \cite{lax}, it is well known 
that through every strictly hyperbolic state $\bar u$ there exist $n$ locally defined and
$C^2$ smooth {\em wave curves}.
These curves collectively provide self-similar solutions to Riemann problems; 
for details see \cite{br}, \cite{daf}, \cite{sm}. Each wave curve is locally made up of 
two components with second order contact at the base point $\bar u$: 
the rarefaction states that are part of the eigencurves, 
and the shock states that are part of the Hugoniot locus $\{\,u\in\Omega\,\,|\,\, \exists\, 
s\in\mathbb R :\, f(u)-f(\bar u)= s(u-\bar u)\,\}$. The geometry of these curves in state 
space thus provide information about the solutions to \eq{claw}.

\subsection{Connections on frame bundles}
\label{geom-prelim}
Given an $n$-dimensional smooth manifold $M$
we let $\mathcal{X}(M)$ and $\mathcal{X}^*(M)$ denote the set of smooth  vector fields and 
differential 1-forms on $M$, respectively.
A \emph{frame} $\{r_1,\dots, r_n\}$  is 
a set of vector fields which span the tangent space $T_p M$ at each point $p\in M$. 
A \emph{coframe} $\{\ell^1,\dots, \ell^n\}$ is a set of $n$ differential 1-forms 
which span the cotangent space $T_p^* M$ at each point $p\in M$.  
The coframe and frame are \emph{dual} if $\ell^i(r_j)=\delta^i_j$ (Kronecker delta). If $u^1,\dots, u^n$ 
are local coordinate functions on $M$, then $\{\pd{}{u^1},\dots,\pd{}{u^n}\}$ is the corresponding local
\emph{coordinate frame}, while $\{du^1,\dots,du^n\}$ is the  dual local \emph{coordinate coframe}.
For a given frame $\{r_1,\dots, r_n\}$ the \emph{structure coefficients} $c^k_{ij}$ are 
defined through
\beq\label{lie-bracket}
	{[r_i,r_j]=\sum_{k=1}^nc^k_{ij}\, r_k}\,,
\eeq
and the dual coframe has related structure equations given by 
\beq\label{dl}
	d\ell^k=-\sum_{i<j} c^k_{ij}\, \ell^i\wedge \ell^j\,.
\eeq 
It can be shown that there exist coordinate functions $w^1,\dots,w^n$ on 
$\Omega$ such that $r_i=\pd{}{w^i},\,i=1\dots,n$, if and only if $r_1,\dots,r_n$ 
commute, i.e.\ all structure coefficients are zero. 
Next, an \emph{affine connection} $\nabla$ on $M$ is an 
$\RR$-bilinear map 
\[\mathcal{X}(M)\times \mathcal{X}(M)\to\mathcal{X}(M)\qquad\qquad
	(X,Y)\mapsto \nabla_X Y\]
such that for any smooth function $f$ on  $M$
\beq\label{connection}
	\nabla_{fX}Y=f\nabla_XY, \qquad \nabla_X (fY)=(Xf)Y+f\nabla_XY\,.
\eeq
By $\RR$-bilinearity and \eq{connection} the connection is uniquely defined by prescribing it 
on a frame:
$$\nabla_{r_i}r_j=\sum_{k=1}^n\Gamma^k_{ij}r_k,$$
where the smooth coefficients $\Gamma^k_{ij}$ are called \emph{connection  components, or Christoffel symbols,} 
relative to the frame $\{r_1,\dots,r_n\}$. Any choice of a frame and $n^3$ functions 
$\Gamma^k_{ij},\, i,j,k=1,\dots,n,$ defines an affine connection on $M$.
A change of frame induces a change of the connection components, and this 
change is not  tensorial. E.g., a connection with zero components 
relative to a coordinate frame, may have non-zero components relative to 
a non-coordinate frame. 

Given a frame $\{r_1,\dots, r_n\}$ with associated Christoffel symbols $\Gamma^k_{ij}$
and dual frame $\{\ell^1,\dots, \ell^n\}$, we define the \emph{connection 1-forms} $\mu_i^j$ by
\beq\label{conn-forms}\mu_i^j:=\sum_{k=1}^n\Gamma^j_{ki}\ell^k\,.\eeq 
In turn, these are used to define two important tensor-fields: the \emph{torsion} 2-forms
\beq\label{tor}
	{\tor}^i:=d\ell^i+\sum_{k=1}^n\mu^i_k\wedge  \ell^k
	=\sum_{k<m}T^i_{km}\ell^k\wedge \ell^m,\qquad i=1,\dots,n\,,
\eeq
and the \emph{curvature} 2-forms
\beq\label{cur}
	\cur^j_i:=d\mu^j_i+\sum_{k=1}^n\mu^j_k\wedge\mu_i^k
	=\sum_{k<m}R^j_{i\,km} \ell^k\wedge \ell^m\,.
\eeq
Here 
\begin{eqnarray}
	\label{T1}T^i_{km}&=&\Gamma^i_{km}-\Gamma^i_{mk}-c^i_{km}\\
	\label{R} R^j_{i\,km}&=&r_k\big(\Gamma^j_{mi}\big)
	-r_m\big(\Gamma^j_{ki}\big)+\sum_{s=1}^n\big(\Gamma^j_{ks}
	\Gamma^s_{mi}-\Gamma^j_{ms}\Gamma^s_{ki}-c^s_{km}
	\Gamma^j_{si}\big)
\end{eqnarray} 
are components of the torsion and curvature  tensors respectively, 
and these \emph{do} change tensorially under a change of frame. 
We can write equations \eq{tor} and \eq{cur} in the compact matrix 
form
\beq\label{tor-cur} 
	\tor=d\ell+\mu\wedge \ell,\quad \cur=d\mu+\mu\wedge\mu
\eeq
where $\ell=(\ell^1,\dots,\ell^n)^T$, $\tor=(\tor^1,\dots,\tor^n)^T$, and 
$\cur$ and $\mu$ are the matrices with components $\cur^j_i$  and $\mu^j_i$ respectively.
The connection is called  \emph{symmetric} if the torsion form is identically zero and  
it is called  \emph{flat}  if the curvature form is  identically zero. Equivalently:
\beq\label{cur0-tor0}	
	d\ell=-\mu\wedge \ell\qquad\mbox{(Symmetry)},
	\qquad
	d\mu=-\mu\wedge\mu \qquad\mbox{(Flatness).}
\eeq
In terms of Christoffel symbols and structure coefficients this 
is equivalent to 
\beq\label{T0}
	c^i_{km} = \Gamma^i_{km}-\Gamma^i_{mk}\qquad\mbox{(Symmetry)}
\eeq
and
\beq\label{R0}
	r_m\big(\Gamma^j_{ki}\big)-r_k\big(\Gamma^j_{mi}\big)
	=\sum_{s=1}^n \big(\Gamma^j_{ks}\Gamma^s_{mi}-\Gamma^j_{ms}
	\Gamma^s_{ki}-c^s_{km}\Gamma^j_{si}\big)\qquad\mbox{(Flatness)}.
\eeq 
One can also show that a connection $\nabla$ is symmetric and flat if and only if in 
a neighborhood  of each point there exist coordinate functions
 $u^1,\dots,u^n$ with the property that the 
Christoffel symbols relative to the coordinate frame are zero: 
$\nabla_{\pd{}{u^i}}\pd{}{u^j}=0$ for all $i,j=1,\dots,n$.

%
%
\subsection{Formulation of the problem}\label{statement}
Returning to \eq{claw} we next provide a precise statement of the ``inverse" problem 
of constructing flux functions $f$ whose 
geometric properties are given. There are various ways to formulate such problems. 
One could prescribe families of curves which are then required to be the Hugoniot loci, the 
eigencurves, or the wave curves for a system of conservation laws \eq{claw}. One might also 
consider giving combinations of these. The most direct formulation is obtained by
prescribing the eigenfields (equivalently, the eigencurves), and this is what we do here.

We will be working locally near a fixed base point $\bar u$ in an open set 
$\Omega\subset \mathbb R^n$, with $\Omega$ smoothly contractible to a point.
Throughout $u=(u^1,\dots,u^n)$ will denote a fixed system of coordinates on a neighborhood 
of $\bar u$. We assume that we are given $n$ linearly independent column $n$-vectors $R_i(u)$ 
(the eigenfields to be) on $\Omega$, and we define
\beq\label{right_left}
	R(u):=[R_1(u)\,|\,\cdots\,|\, R_n(u)]=(R^i_j(u))_{i,\, j}\,,\qquad
	L(u):=R(u)^{-1}=:\left[\begin{array}{c} L^1(u)\\ \hline \vdots\\\hline L^n(u)\end{array}\right]
	=(L^i_j(u))_{i,\, j}\,.
\eeq
The problem we consider may be formulated as follows:
\begin{problem}\label{prob1}
	Consider an open subset $\Omega\subset \mathbb R^n$ equipped with a coordinate system 
	$u=(u^1,\dots,u^n)$. Given a point $\bar u\in \Omega$ together with a (sufficiently smooth) frame 
	$\{R_1(u),\dots,R_n(u)\}$ on $\Omega$. 
	
	Then: find $n$ real functions $\lam^1(u),\dots,\lam^n(u)$ defined on a neighborhood 
	$\mathcal U\subset \Omega$ of $\bar u$ such that, with 
	$\Lambda(u):=\mathrm{diag}[\lam^1(u)\dots\lam^n(u)]$, the matrix
	\beq\label{A}
		A(u):= R(u)\Lambda(u)L(u)
	\eeq
	is the Jacobian matrix with respect to $u$ of some map $f:\mathcal U \to\mathbb R^n$. 
	We are further interested in how large the set of solutions is, i.e.\ how many arbitrary 
	constants and functions appear in a general solution $\lam^1(u),\dots,\lam^n(u)$.
\end{problem}
A solution $\lam(u)=(\lam^1(u),\dots,\lam^n(u))$ to Problem \ref{prob1} generates a
flux $f$ and an associated system of conservation laws \eq{claw} in which 
$u^1,\dots,u^n$ are the {\em conserved quantities}. We note that we do not impose 
strict hyperbolicity for solutions to our problem; indeed, we shall see that there are 
cases where two or more eigenvalues $\lam^i$ must necessarily coalesce.

Below we derive a system of partial differential equations (PDEs) for the 
eigenvalues $\lam^i$, where the coefficients are expressed in terms of the
components of $R_i$ and $L^i$. The system may be formulated in different ways, 
all of which we refer to as the {\em $\lam$-system}. 
Typically this will be an overdetermined system of linear, variable coefficients PDEs.  
It will turn out that the most useful formulation of the $\lam$-system is as an 
{\em algebraic-differential} system; see Section \ref{l_syst}.

\begin{remark}\label{main_rmk}
Let us clarify the coordinate dependence in Problem \ref{prob1}. 
The property of being a Jacobian with respect to a system of coordinates is not invariant 
under most changes of coordinates. For this reason we need to fix the coordinates 
$u^i$ at the outset. 
Of course, having started with a particular choice of coordinates, it may be that the resulting 
$\lambda$-system takes a simpler form when formulated in other coordinates. 
This is the case for so-called {\em rich} systems where the ``good" coordinates are the
Riemann invariants (see Definition \ref{rich_sys} below and Section \ref{rich}). 

It may also be that a change of coordinates happens to preserve Jacobians. This occurs for the 
Euler system of compressible flow when written in Lagrangian variables (Example \ref{gas_dyn}). 
In this case the nonlinear change of dependent variables 
(specific volume, velocity, total energy) $\mapsto$ (specific volume, velocity, entropy)
preserves the conservative form of the system. As far as Problem \ref{prob1} is concerned the two 
forms are equivalent. 
\end{remark}

\begin{remark}
Besides its intrinsic interest Problem \ref{prob1} is motivated by possible applications 
where one searches for systems \eq{claw} with particular geometric properties, and associated properties
(e.g.\ blowup) of solutions to Cauchy problems. In our setup we thus consider the ``geometry" 
(the eigenfields) as given, and we look for the scalar eigenvalue fields.
Alternatively one could prescribe the eigenvalues and search for eigenfields 
$R_i(u)$ that renders $A(u)$ in \eq{A} a Jacobian matrix.
\end{remark}

A special class of systems is given by so-called {\em rich} systems. These are systems 
equipped with a coordinate system of {\em Riemann invariants}. For definitions, and the 
fact that richness can be expressed in terms of the eigenfields, we refer to \cite{serre2}, and 
Section 7.3 in \cite{daf}. Starting from a given frame it is convenient to make the following, equivalent, 
definition:  
\begin{definition}\label{rich_sys}
	The frame $\{R_1(u),\dots,R_n(u)\}$ of linearly independent vector fields is
	said to be {\em rich} if each pair of vector fields is in involution:
	\beq\label{rich_def}
		[R_i,R_j]\in \spa\{R_i,\, R_j\}\qquad \mbox{for all $\, 1\leq i,\, j\leq n$}\,.
	\eeq
\end{definition}
While it is natural to consider the class of rich systems in connection with Problem \ref{prob1}
we shall see that richness does not imply any uniformity with respect to how many constants and
functions are needed to specify a general solution.   

\begin{remark} The class of {\em symmetrizable} systems plays a central role in the
general theory of conservation laws, \cite{fl}, \cite{go}. These are the systems
\eq{claw} that admit a convex ``entropy" $\eta:\Omega\to\RR$, with an associated entropy-flux 
$q:\Omega\to\RR$, such that an additional conservation law
\[\eta(u)_t+q(u)_x=0\]
holds whenever $u$ is a smooth solution of \eq{claw}. Equivalently (at least in the strictly hyperbolic regime),
the Hessian $D^2\eta$ should define an inner product with respect to which the given
eigenvectors $R_i$ are orthogonal ($R_j^T\, D^2\eta\, R_i =0$ for $i\neq j$; see \cite{daf} for details).
Differently from richness, symmetrizability is not expressed directly in terms 
of the eigenfields alone. Problem \ref{prob1} for symmetrizable systems will be taken up 
elsewhere.
\end{remark}

%
%
\subsection{Related works}\label{related}
Problem \ref{prob1} was addressed by Dafermos \cite{daf1} for $2\times 2$-systems 
in several space dimensions under the requirement that the Jacobians
in the various spatial directions commute. Commutativity implies that the Jacobians 
have the same eigenfields. 
In \cite{daf1} it was shown how to construct such systems for any pair of  
linearly independent vector fields. The case of one-dimensional 
$2\times 2$-systems is covered by the analysis of rich systems in Section \ref{rich}.

In his geometric analysis of systems of conservation laws S\'evennec \cite{sev} 
provides a characterization of those quasilinear systems
\[v_t+A(v)v_x=0\]
that can be transformed to conservative form \eq{claw} by a suitable change of variables 
$u=\phi(v)$. In particular, the characterization involves a version of what we refer
to as the $\lam$-system. 

The class of rich systems has been studied by many authors. 
Conlon and Liu \cite{cl} considered rich systems in connection with entropy criteria and 
showed that such systems are endowed with large families of entropies. From a different 
perspective the same class of systems were studied by Tsarev \cites{ts1, ts2}.
S\'evennec \cite{sev} showed that the eigenvalues of strictly hyperbolic, rich systems 
must satisfy certain restrictive conditions (see Proposition 5 in \cite{sev}). Serre \cite{serre2} has 
performed a comprehensive analysis of rich systems, including building of entropies,
commuting families of systems, and construction of rich systems.
In Section \ref{rich-rank0} we analyze rich systems in relation to Problem \ref{prob1}.

%
%
\subsection{Outline and summary of results}\label{outline}
It turns out that the complete solution of Problem \ref{prob1} for arbitrary $n$ is 
quite complicated. In this paper we provide a complete solution for the case $n=3$,
as well as for rich systems for any $n$. (This latter class covers the case $n=2$.)  
We also give a list of examples that illustrate the various types of solutions. 

In Section \ref{l_syst} we begin by noting some properties the $\lam$-system 
that follow directly from the formulation of Problem \ref{prob1}. We then formulate three 
equivalent versions of the $\lam$-system, including an algebraic-differential system, and we record 
the extreme cases with minimal and maximal number of algebraic constraints.
A few general facts are collected in Proposition \ref{gen_case}. 

Our solution of the $n=3$ case in Section \ref{3x3_&_constraints} reveals that 
the solution set of the $\lambda$-system depends on the number of independent algebraic 
equations, as well as on the number of $\lambda^i$ that appear in these equations. 
For $n\geq 4$ there seems to be no easy way to analyze completely the resulting cases.
We provide a complete breakdown of the possible scenarios when $n=3$. 
The algebraic part of the $\lambda$-system now contains two, one, or no independent 
algebraic equations. The two extreme cases fall into either the trivial or rich categories.
In the case of one algebraic equation the size of the solution set depends further on 
the number of $\lambda^i$ involved in the algebraic constraint. The Frobenius integrability
theorem can be applied when all three $\lambda^i$ appear in the algebraic equation, 
and in this case the general solution depends on two constants. 
The only other possibility is that exactly two $\lambda^i$ occur in the algebraic equation. 
This case may be analyzed by using the Cartan-K\"ahler integrability theorem, and the general solution 
now depends on one function of one variable and one constant.

We solve Problem \ref{prob1} for the class of rich systems of any dimension in Section \ref{rich}.
The subclass of rich systems without algebraic constraints are analyzed in Section \ref{rich-rank0}, 
while Section \ref{rich-rank-non-zero} treats the more involved case of rich systems where
the eigenvalues are related algebraically. In both cases the $\lambda$-system can be analyzed
 using an integrability theorem of Darboux (Theorem \ref{dar3}). 

For completeness we include statements of the various integrability theorems we apply.
Concerning smoothness assumptions of the given frame we need to require analyticity 
when we apply the Cartan-K\"ahler theorem. On the other hand, for the cases that 
use the theorems of Frobenius and Darboux we only need to require $C^2$ smoothness 
of the given frame.

Section \ref{examples} collects several examples that are of interest in themselves 
or that illustrate the different cases treated in Section \ref{3x3_&_constraints} and Section \ref{rich}. 
We start by considering the various  
solutions of Problem \ref{prob1} for the case where the given eigenfields are
those of the Euler system describing one-dimensional compressible fluid flow (Example \ref{gas_dyn}).
Depending on the prescribed pressure function these eigenfields may form either a rich or a non-rich frame.
We proceed with several more examples of non-rich frames on $\RR^3$ illustrating  various scenarios 
treated in Section \ref{3x3_&_constraints} and an example of a frame on $\RR^4$ with only trivial 
solutions for the $\lam$-system. We then give a set of examples of  rich frames whose $\lam$-systems 
do not impose algebraic  constraints on eigenfunctions. This includes  a frame on $\RR^2$,  a rich 
orthogonal frame  on $\RR^3$ and a constant frame on $\RR^n$. We finally give two examples of rich 
frames on $\RR^3$ whose  $\lam$-systems impose certain algebraic constraints on eigenfunctions. 
We have written a set of procedures\footnote{ {\sc Maple}  code is available at 
{\tt http://www.math.ncsu.edu/$\sim$iakogan/mapleHTML/lambda-system.html}} in {\sc Maple} 
to obtain and analyze the $\lam$-system for a given frame, and used it
to construct some of the above examples.
 
Finally, our analysis indicates that a solution of Problem \ref{prob1} for general systems with $n\geq 4$ 
is rather involved, with a large number of different subcases.

%
%
\section{The $\lam$-system}\label{l_syst}
\subsection{Trivial solutions and scalings}\label{triv_scal}
Problem \ref{prob1} always has a one-parameter family of solutions 
given by
\[\lam^1(u)=\cdots=\lam^n(u)\equiv \hat \lam\,,\]
where $\hat\lam$ is any real constant. We refer to these as {\em trivial} solutions.
The resulting system \eq{claw} is linearly degenerate in all families (i.e.\ 
$\nabla \lam_i(u)\cdot R_i(u)\equiv 0$) and  any map $f(u)=\hat\lam u +\hat u$, 
where $\hat u\in\mathbb R^n$, is a corresponding flux. 
While not of interest themselves, the trivial solutions
show that any compatibility condition associated with the $\lam$-system will
not rule out existence of solutions altogether.
A trivial solution can be added to any solution of the $\lam$-system 
to give another solution of the same system.  
For later reference we record the following related result:
\begin{proposition}\label{equal_lambdas}
	If $\lam^1(u)=\cdots=\lam^n(u)$ is a solution
	to Problem \ref{prob1}, then their common value is a constant.
\end{proposition}
\begin{proof}
	If $\lam^1(u)=\cdots=\lam^n(u)=\hat\lam(u)$ then 
	$A(u)=\Lambda(u)=\hat\lam(u)I_{n\times n}$.
	For this to be a Jacobian of a map $f:\mathcal U \to\mathbb R^n$ 
	we must have that
	\[\forall\, i,j=1,\dots,n\,:\qquad  \frac{\del f_i}{\del u^j}(u)
	=\delta_{ij}\hat\lam (u)\,,\]
	whence $\hat\lam(u)$ is a function of $u^i$ alone for each $i$, 
	and thus constant.
\end{proof} 
In formulating Problem \ref{prob1} we may 
use any (non-vanishing) scalings of the eigenfields $R_i(u)$.  
That is, given smooth functions
$\alpha^j:\Omega\to \mathbb R\setminus \{0\}$, $j=1,\dots, n$, we may set 
$\ti R_j(u)= \alpha^j(u)R_j(u)$, together with the inversely scaled left eigenfields 
$\ti L_j(u):= {\alpha^j(u)}^{-1}L_j(u)$. Letting $\ti R:=R\alpha$, $\ti L:=\alpha^{-1}L$, 
where $\alpha(u)=\mbox{diag}[\alpha_1(u)\dots\alpha_n(u)]$, we get that  
$R(u)\Lambda(u) L(u)$ is a Jacobian if and 
only if $\ti R(u)\Lambda(u) \ti L(u)$ is a Jacobian. 
E.g., in the case of rich systems 
the simplest form of the $\lam$-system is obtained by using the versions 
that makes the matrix $L$ a Jacobian with respect to the $u$ variable 
(see Section \ref{rich}).

\subsection{Formulating the $\lam$-system}
\subsubsection{Direct formulation}\label{direct}
We have that a matrix $A(u)=(A^i_j(u))_{i,\, j}$, 
is a Jacobian with respect to $u$-coordinates if and only if 
\beq\label{sys1}
    \del_k A^i_j(u) = \del_j A^i_k(u)
    \qquad\qquad\mbox{for all $\, i, j, k=1,\dots,n\, $ with $\, j<k$}\,,
\eeq
where $\del_i$ denotes partial differentiation with respect to $u_i$.
In the case when $A(u)$ is given by \eq{A} we set
\[C^i_{mj}(u):=R^i_m(u)L^m_j(u)\qquad\mbox{(no summation),}\]
and \eq{sys1} may be written
\beq\label{sys2}
    \sum_{m=1}^n \Big[ C^i_{mj}\del_k \lam^m - C^i_{mk}\del_j \lam^m
    + \lam^m\big(\del_k C^i_{mj} - \del_j C^i_{mk} \big)\Big] = 0\,,
\eeq
where $i,j,k\in\{1,\dots, n\} $ with $j<k$.
This version of the $\lam$-system provides a homogeneous, variable coefficient
system of $\frac{n^2(n-1)}{2}$ linear  PDEs for $n$ unknowns.
For $n\geq 3$ it is thus typically an overdetermined system of PDEs.

\subsubsection{Formulations using 1-forms}\label{forms}
A simpler formulation of the $\lam$-system is obtained by using differential 
forms to express the condition that $A(u)$ in Problem \ref{prob1} is a 
Jacobian. By Poincar\'e's lemma (recall that $\Omega$ is assumed smoothly 
contractible to a point) we have
\beq\label{cond}
    \mbox{$A(u)$ is a Jacobian matrix w.r.t.\ $u$}
    \qquad \Longleftrightarrow\qquad dA(u)\wedge du = 0\,,
\eeq
where the $d$-operator is applied component-wise.
Applying the product rule, condition \eq{cond} is thus 
equivalent to
\begin{equation}\label{cond_I}
    \left\{L(dR)\Lambda + d\Lambda -
    \Lambda L(dR) \right\}\wedge Ldu=0\,,
\end{equation}
or
\begin{equation}\label{cond_I_b}
    	(d\Lambda)\wedge (Ldu) = 
	\big\{\Lambda (L dR)-(L dR)\Lambda \big\}\wedge (Ldu)\,,
\end{equation}
where we have used that $L=R^{-1}$. (Clearly \eq{cond_I_b} is satisfied 
if $d\Lambda=\Lambda (L dR)-(L dR)\Lambda$; however, the associated 
solutions are exactly the trivial solutions $(\lam_1,\dots, \lam_n)\equiv 
(\hat\lam,\dots, \hat\lam)$, $\hat \lam\in\RR$.)

The system \eq{cond_I_b} is an equation for $n$-vectors of 2-forms.  
We proceed to write out the system in $u$-coordinates by applying \eq{cond_I_b} 
to the pair of vector fields $(\frac{\del}{\del u^i},\frac{\del}{\del u^j})$.  A direct 
calculation yields:
\beq\label{cond_III}
	L^i_j(\del_l\lam^i) - L^i_l(\del_j\lam^i)=
	\sum_{m\neq i}B^{mi}_{lj} (\lam^m-\lam^i)\,,\qquad
	\mbox{for each $i$ and for all $1\leq l < j \leq n$,}
\eeq
where $B^{mi}_{lj}:= \big\{L^m_j (\del_l L^i) - L^m_l(\del_j L^i)\big\}\cdot R_m$.
Not surprisingly this is again a homogeneous system of $\frac{n^2(n-1)}{2}$ 
linear PDEs for $n$ unknowns.
Thus, in this formulation of the $\lam$-system the 
equations involving derivatives of $\lam^i$ are such that, first, 
no other eigenvalue appears differentiated, and second, at most two 
derivatives of $\lam^i$ occur in each equation.

To obtain a formulation of the $\lam$-system which explicitly 
records {\em algebraic} constraints on the eigenvalues $\lam^i$,
we introduce the following notation. Whereas we treat $R_i(u)$ 
and $L^i(u)$ as arrays of functions we let $r_i(u)$, $\ell^i(u)$ denote 
the corresponding vector fields (differential operators) and differential 
1-forms:
\beq\label{v_flds,1-frms}
	r_i(u):=\sum_{m=1}^n R_i^m(u)\frac{\del{\quad\! }}{\del u^m}\Big|_u\,,
	\qquad \ell^i(u):=\sum_{m=1}^n L_m^i(u)du^m\big|_u\, .
\eeq
Since $R(u)$ is assumed invertible on $\Omega$, the vector fields 
$r_i(u),\,i=1,\dots,n$, provide a frame on $\Omega$, and the 1-forms 
$\ell^i(u),i=1,\dots,n$, provide the dual coframe on $\Omega$ 
(see Section~\ref{geom-prelim}).  We define the $n$-vector of 1-forms $\ell$ by
\[\ell:=\left[\begin{array}{c}
	\ell^1\\
	\vdots\\
	\ell^n
\end{array}\right]=Ldu\,,\]
and  introduce the following coefficients
\beq\label{christoffel}
	\Gamma_{ij}^k:=L^k(D R_j)R_i\,,
\eeq
where $D$ denotes Jacobian with respect to $u$. 
A direct computation shows that they are, in fact, the Christoffel symbols 
(connection components) of the flat and symmetric connection 
$\nabla_{\frac{\del{\,}}{\del u^i}}\textstyle\frac{\del{\,}}{\del u^j}=0\, $
computed  relative to the frame $\{r_1,\dots,r_n\}$. In other words, the 
covariant derivatives of the frame vector fields are 
$\nabla_{r_i}{r_j}=\sum_{k=1}^n\Gamma_{ij}^k r_k$.
We note that  the symmetry condition \eq{T0} implies
\beq\label{r_comm}
	[r_i,r_j]=\sum_{k=1}^nc^k_{ij}\, r_k= \sum_{k=1}^n(\Gamma_{ij}^k-\Gamma_{ji}^k)r_k\,.
\eeq
This shows that if $\{r_1,\dots,r_n\}$ is rich, then
$\Gamma_{ij}^k=\Gamma_{ji}^k$ whenever $k\notin\{i,\,j\}$.

Another direct  calculation  shows that the $(k,j)$-entry of the 
matrix $\mu:=R^{-1}dR=LdR$ of 1-forms is given by the connection forms
\beq\label{mu}
	\mu^k_j:=\big(LdR\big)_j^k=\sum_{i=1}^n \Gamma_{ij}^k \ell^i\,.
\eeq
Thus \eq{cond_I_b} reads
\beq\label{lambda-in-r-frame}
	{(d\Lambda)\wedge \ell = 
	\big\{\Lambda \mu-\mu\Lambda \big\}\wedge \ell\,}.
\eeq
Again \eq{lambda-in-r-frame} is an equation of $n$-vectors of 2-forms. Applying
each component to pairs of frame vector fields $(r_i,r_j),\, i,j=1,\dots ,n$, 
we obtain an equivalent form of the $\lambda$-system as a differential-algebraic system:
\begin{eqnarray}
	r_i(\lam^j) &=& \Gamma_{ji}^j (\lam^i-\lam^j)
	\qquad \mbox{for $i\neq j$,}\label{sev1}\\
	(\lam^i-\lam^k)\Gamma_{ji}^k &=& 
	(\lam^j-\lam^k)\Gamma_{ij}^k
	\qquad \mbox{for $i<j,\, i\neq k,\, j\neq k$,}\label{sev2}
\end{eqnarray}
where there are no summations. \eq{sev1} gives $n(n-1)$ linear, homogeneous
PDEs, while \eq{sev2} gives $\frac{n(n-1)(n-2)}{2}$ algebraic relations. We observe that
\eq{sev2} is symmetric in $i$ and $j$, and that {\em all} coefficients 
$\Gamma_{ij}^k$ with $i\neq j\neq k\neq i$ appear in \eq{sev2}.
This form of the $\lam$-system 
explicitly records the algebraic relations inherent in \eq{sys2} and \eq{cond_III}. 
\begin{remark}\label{sev_rmk}
	The equations \eq{sev1}-\eq{sev2} appear in S\'evennec's 
	characterization of quasilinear systems that admit a conservative 
	form, see  \cite{sev}. For a given quasilinear system
	\[v_t+A(v)v_x=0\,,\qquad A(v)\in \RR^{n\times n}\,,\]
	S\'evennec shows that there is a coordinate system in which the 
	system is conservative if and only if there exists a flat and symmetric affine connection 
	$\nabla$ such that its Christoffel symbols and the eigenvalues of $A(u)$ satisfy 
	\eq{sev1}-\eq{sev2}. 
	The same system appears in Tsarev \cite{ts2}.
\end{remark}
\begin{remark}\label{diff-closed}
Note that condition \eq{lambda-in-r-frame} is equivalent to
\beq\label{Theta}
	\Theta:=\big(d\Lambda -
	\Lambda \mu+\mu\Lambda \big)\wedge \ell\,=0\,,
\eeq
where $\Theta$ is a differential form on a $2n$-dimensional manifold with 
coordinates $u^1,\dots,u^n,\lambda^1,\dots, \lambda^n$. Moreover, the ideal $\mathcal{I}$ algebraically 
generated by $\Theta$ in the ring of
differential forms  is a differential ideal, i.e. $d\omega\in \mathcal{I}$ for all $\omega\in \mathcal{I}$.
Indeed, from the flatness and symmetry properties of the connection:
\[d\mu=-\mu\wedge\mu \qquad\mbox{and}\qquad dl=-\mu\wedge l\,,\] 
it follows that
\begin{eqnarray*}d\Theta&=&\big(-d\Lambda\wedge \mu-\Lambda(d\mu)+
	(d\mu)\Lambda-\mu\wedge d\Lambda \big)\wedge\ell - \big(d\Lambda -
	\Lambda \mu+\mu\Lambda \big)\wedge d\ell\\
	&=&\big(-d\Lambda\wedge \mu+\Lambda(\mu\wedge\mu)-(\mu\wedge\mu)
	\Lambda-\mu\wedge d\Lambda \big)\wedge\ell + \big(d\Lambda -
	\Lambda \mu+\mu\Lambda \big)\wedge\mu\wedge\ell\\
	&=& \big(-(\mu\wedge\mu)\Lambda-\mu\wedge d\Lambda\big)\wedge\ell+(\mu\Lambda)
	\wedge\mu\wedge \ell\\
	&=& \mu\wedge\big(-\mu\Lambda- d\Lambda+\Lambda\,\mu\big)\wedge\ell
	=-\mu\wedge\Theta\equiv 0 \mbox{ mod } \Theta\,.
\end{eqnarray*}
Solving the $\lam$-systems is therefore equivalent to finding $n$-dimensional integral 
sub-manifolds $\lambda^i=F^i(u^1,\dots,u^n),\,i=1,\dots,n$, of $\mathcal I$. 
\end{remark}

\subsection{Observations about the rank of the algebraic sub-system \eq{sev2}}\label{rank_discuss}
We next consider the algebraic constraints \eq{sev2}, which is a system of $\frac{n(n-1)(n-2)}{2}$ 
linear equations. Choosing the variables to be the differences
\[x^k:=\lam^k-\lam^1\,,\qquad \qquad k=2,\dots,n\,,\]
we rewrite \eq{sev2} in matrix form as
\beq\label{x_syst}
	Nx=0\,.
\eeq
Here $x$ is  the $(n-1)$-vector $(x^2,\dots,x^n)^T$ and $N$ is a certain  
$\frac{n(n-1)(n-2)}{2}\times (n-1)$-matrix whose entries are given in terms of 
the $\Gamma_{ij}^k$. It is easily checked that each entry of $N$ is either zero, a single Christoffel 
symbol ($\pm\Gamma_{ij}^k$), or a difference of such ($\Gamma_{ij}^k-\Gamma_{ji}^k$).
Furthermore, $N=0$ if and only if $\Gamma_{ij}^k=0$ for all choices $i\neq j\neq k\neq i$.

The number of independent algebraic constraints is given by $\rank(N)$.
It is convenient to use this as a first, rough classification. However, $\rank(N)$ does not characterize 
the solutions to the $\lam$-system (in terms of number of constants and functions in the general 
solution). This is clear already from the cases of maximal and minimal rank. 
We briefly consider these extreme cases:
\begin{itemize}
	\item $\rank(N)=0$. In this case $N=0$ and there are no algebraic constraints 
	imposed on the eigenvalues. This occurs if and only if $\Gamma_{ij}^k=0$ for all choices of 
	$i\neq j\neq k\neq i$. By \eq{r_comm}
	this implies that $[r_i,r_j]\in\spa\{r_i,r_j\}$, i.e.\ we are in the rich case: 
	\[\rank(N)=0\quad \Rightarrow \quad \{r_1,\dots,r_n\}\mbox{\, is rich}\,.\] 
	The case 
	$\rank(N)=0$ will be analyzed in more detail in Section \ref{rich-rank0}. 
	We will show that any frame on $\RR^2$, any constant frame on $\RR^n$, as well as any rich 
	orthogonal frame on $\RR^n$, falls in this category. Corresponding
	examples are given in Section~\ref{ex_rich}. 
	  
	The question of whether richness implies $\rank(N)=0$ is somewhat subtle. 
	Namely, if we were given a rich and strictly hyperbolic system \eq{claw},
	then indeed $\rank(N)=0$. (See \cite{daf} p.\ 185 for the proof that $\Gamma_{ij}^k=L^k(D R_j)R_i=0$
	for all choices of $i\neq j\neq k\neq i$ in this case.)
	However, with Problem \ref{prob1} we are starting from given vector fields $r_i$ without 
	insisting on strict hypebolicity.
	Example~\ref{rich_rank1} and Example~\ref{rich_rank2} show that it is possible to prescribe a collection of 
	vector fields that form a rich family, without $\rank(N)$ being zero. 
	Furthermore,  Example~\ref{rich_rank1} shows that   the associated 
	$\lam$-system for frames of this type may admit nontrivial solutions.
	Thus:
	\[\{r_1,\dots,r_n\}\mbox{\, is rich}\quad \not\Rightarrow \quad \rank(N)=0\,.\]
	 In Section~\ref{rich-rank-non-zero} we prove, however, that
	the
	$\lam$-system  associated with a rich, $\rank(N)>0$ frame allows no strictly hyperbolic solutions.
	 \item $\rank(N)=n-1$. In this case the only solution to  
	\eq{x_syst} is $x=0$, that is, all $\lam^i$ are equal. According to Proposition
	\ref{equal_lambdas} it follows that all eigenvalues are equal to a 
	common {\em constant}:
	$$\rank(N)=n-1\quad \Rightarrow \quad\mbox{$\lam$-system has only trivial solutions}.$$
	 In particular, if the $\lam$-system admits a strictly 
	hyperbolic solution, then necessarily 
	$\rank(N)<n-1$. It is a non-obvious fact that there {\em are} cases where all
	solutions to \eq{sev1}-\eq{sev2} are trivial; explicit examples are provided by 
	Examples ~\ref{nonrich_IIa_trivial}, ~\ref{nonrich_IIb_trivial}, ~\ref{nonrich_n4_maximalrank} and ~\ref{rich_rank2}.
	
	On the other hand we observe that the condition $\rank(N)=n-1$ does not characterize 
	the cases where the  $\lam$-system \eq{sev1}-\eq{sev2} has only trivial solutions. 
	In other words, it may be that the only solutions to \eq{sev1}-\eq{sev2} are the trivial 
	ones, while $\rank(N)<n-1$; for a concrete example see Example~\ref{nonrich_IIa_trivial} 
	and ~\ref{nonrich_IIb_trivial}. Thus:  
	$$\mbox{$\lam$-system has only trivial solutions} \quad \not\Rightarrow \quad \rank(N)=n-1.$$	
\end{itemize}
We summarize our findings:
\begin{proposition}\label{gen_case}
	Consider Problem \ref{prob1} for a given frame $\{R_1(u),\dots,R_n(u)\}$. Then the
	$\lam$-system may be formulated as 
	an algebraic-differential system \eq{sev1}-\eq{sev2} for the eigenvalues $\lam^i$. 
	Absence of the algebraic constraints \eq{sev2} implies that the frame  is rich, 
	but not vice versa. A maximal number of $n-1$ independent algebraic constraints 
	\eq{sev2} implies that the $\lam$-system admits only trivial solutions, but not vice 
	versa.
\end{proposition}

%
%
\section{Systems of three equations}
\label{3x3_&_constraints}
In this section we present a complete breakdown of the possible
solutions of Problem \ref{prob1} for a  given frame of vector-fields in $\RR^3$. 
When $n=3$ the $\lambda$-system \eq{sev1}-\eq{sev2} consists 
of six linear PDEs and three linear algebraic equations. For concreteness 
we record these; the PDEs are
\bea
	r_1(\lam^2) &=& \Gamma_{21}^2(\lam^1-\lam^2) \label{pde1}\\
	r_1(\lam^3) &=& \Gamma_{31}^3(\lam^1-\lam^3)\label{pde2}\\ 
	r_2(\lam^1) &=& \Gamma_{12}^1(\lam^2-\lam^1)\label{pde3}\\
	r_2(\lam^3) &=& \Gamma_{32}^3(\lam^2-\lam^3)\label{pde4}\\ 
	r_3(\lam^1) &=& \Gamma_{13}^1(\lam^3-\lam^1)\label{pde5}\\
	r_3(\lam^2) &=& \Gamma_{23}^2(\lam^3-\lam^2)\,,\label{pde6}
\eea
while the algebraic constraints may be written as
\beq\label{3x3_alg_cond}
	Nx=\left[\begin{array}{cc}
	\Gamma_{32}^1 & -\Gamma_{23}^1\\
	(\Gamma_{31}^2-\Gamma_{13}^2) &  \Gamma_{13}^2\\
	\Gamma_{12}^3 & (\Gamma_{21}^3-\Gamma_{12}^3)
	\end{array}\right]
	\left[\begin{array}{c}
	x^2\\
	x^3
	\end{array}\right]=0\,,
\eeq
where $x^2=\lam^2-\lam^1$ and $x^3=\lam^3-\lam^1$.
There are three possibilities depending on $\rank (N)$:
\begin{itemize}
	\item [I:] $\rank(N)=0$. There are no algebraic constraints; the eigenvectors are 
	pairwise in involution, and any corresponding system of conservation laws 
	\eq{claw} is rich. The analysis in Section \ref{rich-rank0}, which applies to systems of any size,
	demonstrates that the $\lam$-system always has many non-trivial (in particular, many 
	strictly hyperbolic) solutions in this case. 
	\item [II:] $\rank(N)=1$. In this case \eq{3x3_alg_cond} imposes a single 
	linear relationship among the eigenvalues. This case is more involved and 
	there are several possibilities in terms of how many constants and functions determine 
	a general solution. The analysis is detailed in Section \ref{rank1_n=3} below.
	\item [III:] $\rank(N)=2$. In this case there are only trivial solutions 
	$\lam^1=\lam^2=\lam^3\equiv$ constant. 
\end{itemize}
Section \ref{examples} provides examples for each type of behavior.

\subsection{Case II: a single algebraic relation}\label{rank1_n=3}
Using \eq{T0} we rewrite the algebraic relations \eq{3x3_alg_cond}:
\bea
	c_{32}^1\lam^1&=& \Gamma_{32}^1\lam^2-\Gamma_{23}^1\lam^3\,, \label{c1}\\
	c_{31}^2\lam^2&=& \Gamma_{31}^2\lam^1-\Gamma_{13}^2\lam^3\,, \label{c2}\\
	c_{21}^3\lam^3&=& \Gamma_{21}^3\lam^1-\Gamma_{12}^3\lam^2\,. \label{c3}
\eea
By assumption the rank of the system \eq{3x3_alg_cond} is 1, whence the three equations are 
all equivalent to a single non-trivial algebraic condition (unique up to non-vanishing scalings)
\beq\label{alpha-lambda}
	\alpha_1\lam^1+\alpha_2\lam^2+\alpha_3\lam^3=0\,,
\eeq
where necessarily
\beq\label{alphas}
	\alpha_3=-(\alpha_1+\alpha_2)\,.
\eeq
There are therefore two sub-cases to consider: 
\begin{itemize}
\item [$\bullet$ {\bf IIa:}] all three $\lambda^i$ appear in \eq{alpha-lambda} 
with non-zero coefficients,
\item [$\bullet$ {\bf IIb:}] only two of three $\lambda^i$ are involved in \eq{alpha-lambda} 
with non-zero coefficients. 
\end{itemize}
In either case it may be that the only solutions are trivial.
To analyze non-trivial solutions we employ the Frobenius integrability theorem 
in case IIa, while case IIb requires the more general Cartan-K\"ahler integrability 
theorem. It will turn out that the number of non-trivial solutions differ in the two situations.

\subsection{Subcase IIa: All three $\lambda^i$ appear 
in the unique algebraic relation}
We first recall the relevant formulation of the Frobenius integrability theorem.
\begin{definition}\label{Frob_def}
Let $\mathcal U\subset\RR^m,\,  \mathcal V\subset \RR^n$ be open sets, and
assume that $g^i_j:\mathcal U\times \mathcal V\to\RR$, $1\leq j\leq m$, $1\leq i\leq n$, are  
smooth functions and $\{Y_1,\dots,Y_m\}$ is a frame on $\mathcal U$.
Then the first order PDE system 
\beq\label{g_syst}
	Y_j(v^i)=g^i_j(x,v(x))\qquad\qquad\mbox{for the unknown $n$-vector 
	$v(x)=(v^1(x),\dots,v^n(x))$,}
\eeq
is called a {\em Frobenius system} (in $n$ unknowns on $\mathcal U\times \mathcal V$). 
\end{definition}
That is, a Frobenius system prescribes {\em all} first derivatives of {\em all} the unknowns,
see \cite{sp}.

\begin{theorem} [Frobenius Integrability Theorem - Frame Version] Suppose  
\eq{g_syst} is a Frobenius system in $n$ unknowns on $\mathcal U\times \mathcal V$ 
that satisfies the following integrability conditions as identities in $(x,v)\in\mathcal U\times \mathcal V$:
\beq\label{cc}
	\hat Y_k\big(g_j^i\big) - \hat Y_j\big(g_k^i\big) =  \sum_{l=1}^m \alpha_{kj}^l g^i_l
	\qquad\qquad\mbox{for}\quad 1\leq j < k\leq m,\,\,1\leq i\leq n\,.
\eeq
Here the structure coefficients $\alpha_{kj}^l$, $1\leq j,\, k,\, l\leq m$, are given by
$[Y_k,Y_j]=\sum_{l=1}^m \alpha_{kj}^l Y_l$, and $\hat Y_k\big(g_j^i\big)$ denote the total 
derivatives obtained by using the equations in \eq{g_syst}:
$\hat Y_k\big(g_j^i\big):= Y_k (g_j^i) + \big(\nabla_u g_j^i\big)\cdot g_k\,.$
Then, for a fixed point $(\bar x,\bar v)\in\mathcal U\times \mathcal V$ (and under suitable 
smoothness conditions on $g^i_j$, $Y_k$), the system \eq{g_syst} has a unique local solution
$v(x)$ defined for $x$ near $\bar x$. Furthermore, these solutions foliate a neighborhood of 
$(\bar x,\bar v)$ as $\bar v$ varies over $\mathcal V$, and the general solution to 
\eq {g_syst} depends on $n$ constants.
\end{theorem}

\noindent
In the case IIa it turns out that the $\lambda$-system can be rewritten as a Frobenius system 
and we have:
\begin{theorem}\label{case_IIa}
	Assume $n=3$ and that the $\lambda$-system contains a single algebraic constraint 
	\eq{alpha-lambda} ($\rank(N)=1$). 
	Consider the case IIa where all three $\lambda^i$ appear with non-vanishing coefficients in 
	\eq{alpha-lambda}. Then, after elimination of one of the unknowns, the PDEs in the 
	$\lambda$-system, reduce to a Frobenius system for two unknown functions of three variables.  
	If the corresponding compatibility conditions \eq{cc} are satisfied as identities then the general 
	solution to the $\lambda$-system \eq{pde1}-\eq{pde6}, \eq{3x3_alg_cond} depends on two 
	constants. Otherwise, the only solutions are trivial solutions. Finally, both situations can occur.
\end{theorem}
\begin{proof}
We return to \eq{alpha-lambda}, where we assume that all three $\lambda^i$ are non-vanishing.
By relabeling indices if necessary we may assume that \eq{alpha-lambda} is proportional to  \eq{c1}, 
and hence  $c_{32}^1\neq 0$,  $\Gamma_{32}^1\neq 0$ and $\Gamma_{23}^1\neq 0$. 
We solve \eq{c1} for $\lam^1$:
\beq\label{solve1}
	\lam^1=\frac 1{c_{32}^1}(\Gamma_{32}^1\lam^2-\Gamma_{23}^1\lam^3)\,,
\eeq
and use this to eliminate $\lam_1$ in the differential equations \eq{pde1}-\eq{pde6}. 
By using repeatedly that $c_{ij}^k=\Gamma_{ij}^k-\Gamma_{ji}^k$, such that 
$r_m\big(\frac{\Gamma_{ij}^k}{ c_{ij}^k}\big)=r_m\big(\frac{\Gamma_{ji}^k}{ c_{ij}^k}\big)$ 
for all $i,j,k,m$, we obtain the PDEs
\begin{eqnarray}
\label{e1}	r_1(\lam^2) &=& \frac{\Gamma_{21}^2\Gamma_{23}^1}{c_{32}^1}(\lam^2-\lam^3)\\
\label{e2}	r_1(\lam^3) &=& \frac{\Gamma_{31}^3\Gamma_{32}^1}{c_{32}^1}(\lam^2-\lam^3)\,,\\
\label{e3}	r_2\left(\frac{\Gamma_{32}^1}{c_{32}^1}\right) (\lam^2-\lam^3) 
		+\frac {\Gamma_{32}^1}{c_{32}^1}r_2(\lam^2)-\frac{\Gamma_{23}^1}{c_{32}^1}r_2(\lam^3) 
		&=& \frac{\Gamma_{12}^1\Gamma_{23}^1}{c_{32}^1}(\lam^3-\lam^2)\\
\label{e4}	r_2(\lam^3) &=& \Gamma_{32}^3(\lam^2-\lam^3)\,,\\
\label{e5}	r_3\left(\frac{\Gamma_{32}^1}{c_{32}^1}\right) (\lam^2-\lam^3) 
		+\frac {\Gamma_{32}^1}{c_{32}^1}r_3(\lam^2)-\frac{\Gamma_{23}^1}{c_{32}^1}r_3(\lam^3)  
		&=& \frac{\Gamma_{13}^1\Gamma_{32}^1}{c_{32}^1}(\lam^3-\lam^2)\\
\label{e6}	r_3(\lam^2) &=& \Gamma_{23}^2(\lam^3-\lam^2)\,,
\end{eqnarray}
Since $\Gamma_{32}^1,\, \Gamma_{23}^1\neq 0$, we can solve \eq{e3} and \eq{e5} for $r_2(\lam^2)$ and 
$r_3(\lam^3)$ by using \eq{e4} and \eq{e6}:
\begin{eqnarray}
	\nonumber r_1(\lam^2) = \frac{\Gamma_{21}^2\Gamma_{23}^1}{c_{32}^1}(\lam^2-\lam^3)\,,&&\qquad
	 r_1(\lam^3) = \frac{\Gamma_{31}^3\Gamma_{32}^1}{c_{32}^1}(\lam^2-\lam^3)\,,\\
	\qquad \label{IIa-solved} r_2(\lam^2) = \left[\frac{\Gamma_{23}^1}{\Gamma_{32}^1}(\Gamma_{32}^3-\Gamma_{12}^1)
	-\frac{c_{32}^1}{\Gamma_{32}^1}\,r_2\left(\frac{\Gamma_{32}^1}{c_{32}^1}\right)\right](\lam^2-\lam^3)\,,&&\qquad
	r_2(\lam^3) = \Gamma_{32}^3(\lam^2-\lam^3)\,,\\
	\nonumber  r_3(\lam^3) =  \left[\frac{\Gamma_{32}^1}{\Gamma_{23}^1}(\Gamma_{13}^1-\Gamma_{23}^2)
	+\frac{c_{32}^1}{\Gamma_{23}^1}\,r_3\left(\frac{\Gamma_{23}^1}{c_{32}^1}\right)\right](\lam^2-\lam^3)\,,&&\qquad
	r_3(\lam^2) =- \Gamma_{23}^2(\lam^2-\lam^3)\,.
\end{eqnarray}
This system specifies the derivatives of the two unknown functions $\lam^2$ and $\lam^3$ along
all three vector fields $r_1$, $r_2$ and $r_3$. Hence the system is of Frobenius type. 
For simplicity we write the system as
\beq\label{simple}
	r_i(\lambda^s)=\phi_i^s(u)(\lambda^2-\lambda^3)\qquad \mbox{for $i=1,2,3$ and  $s=2, 3,$}
\eeq
where $\phi_i^s$ are known functions of $\Gamma$'s, given by the right-hand sides in \eq{IIa-solved}. 
According to the Frobenius Integrability Theorem this system is integrable provided 
\beq\label{comp}
	[r_i,r_j]=\sum_{k=1}^3c^k_{ij}\, r_k\,,\qquad \mbox{for $1\leq i<j\leq 3$}\,,
\eeq
where $c_{ij}^k=\Gamma_{ij}^k-\Gamma_{ji}^k$ and the left-hand side is computed 
by using the equations \eq{simple}. A calculation reduces \eq{comp} to:
\beq\label{comp2}
	\Big[r_i(\phi_j^s)-r_j(\phi_i^s)+\phi_j^s(\phi_i^2-\phi_i^3)-\phi_i^s(\phi_j^2-\phi_j^3)\Big](\lam^2-\lam^3)=
	\left[\sum_{k=1}^3c^k_{ij}\phi_k^s\right](\lam^2-\lam^3)\,, 
\eeq
where $1\leq i<j\leq 3$ and $s=2,3$.

These conditions are satisfied if $\lam^2=\lam^3$, in which case the system \eq{IIa-solved} implies that 
$\lam^2=\lam^3$ is a constant. Equation \eq{solve1} then shows that 
$\lam^1=\lam^2$, and Proposition \ref{equal_lambdas} implies that the solution in this case is trivial: 
$\lam^1=\lam^2=\lam^3\equiv constant$.

For a non-trivial  solution  to exist the following six conditions must hold:
\begin{eqnarray}
\label{frob-comp1}	r_i(\phi_j^2)-r_j(\phi_i^2)
	&=&\phi_j^2\phi_i^3-\phi_i^2\phi_j^3+\sum_{k=1}^3c^k_{ij}\phi_k^2\qquad 1\leq i<j\leq 3,\label{s=2} \\
	\label{frob-comp2}r_i(\phi_j^3)-r_j(\phi_i^3)&=&
	\phi_j^2\phi_i^3-\phi_i^2\phi_j^3+\sum_{k=1}^3c^k_{ij}\phi_k^3\qquad 1\leq i<j\leq 3.\label{s=3} 
 \end{eqnarray}
Example~\ref{gas_dyn} and Example~\ref{nonrich_IIa_trivial} show that these compatibility conditions may or may 
not be satisfied: they must be checked for each case individually. 
If the compatibility conditions are met then, according to the Frobenius Theorem, the general solution
to the $\lambda$-system depends on two constants. 
\end{proof}

\subsection{Subcase IIb: Exactly two $\lambda^i$ appear 
in the unique algebraic relation}\label{subcaseIIb}
This case is more involved than IIa: the $\lambda$-system does not reduce to a 
Frobenius system and the non-trivial solutions must be analyzed by using the more 
general Cartan-K\"ahler theorem. This requires a reformulation of the $\lam$-system as an 
Exterior Differential System (EDS). Using the terminology and notation of  \cite{bcggg} and \cite{il}
we give the following  formulation of  Cartan-K\"ahler theorem (see Theorem 7.3.3  in \cite{il} and 
discussion on p. 87 in \cite{bcggg}). In order to apply this theorem we need to assume that the 
given frame is analytic.

\begin{theorem}(Cartan-K\"ahler Integrability Theorem). \label{Cartan-Kahler} Let $E_{k},\, 0\leq k\leq n$, 
be a flag of integral elements at a point $p$ for an analytic EDS  on an $(n+s)$-dimensional manifold, 
with $\dim E_k=k$, and such that  $E_k$ is K\"ahler regular. Then there exists a smooth $n$-dimensional 
integral manifold $S$ whose tangent space at $p$ is $E_n$.
Furthermore, let $H(E_k)$ be the polar spaces of $E_k$ and  $c_k=\codim H(E_k)$ for $0\leq k \leq n-1$, 
with $c_{n}=\codim E_{n}$. Let $x=(x^{1},\dots,x^{n})$ and $y=(y^{1},\dots, y^{s})$ be local coordinates around  
$p$, chosen so that  $E_{k}$ is spanned by $\dbyd{x^{1}},\dots,\dbyd{x^{k}}$, $E_{n}$ is annihilated by  
$dy^{1}, \dots,dy^{n}$, and $H(E_{k})$ is annihilated by 
$dy^{1}, \dots,dy^{c_{k}}$, for  $0\leq k \leq n-1$.  
Then $S$ is defined in the coordinates $(x,y)$ by analytic equations of the form
$y^{\alpha}=F^{\alpha}(x^{1},\dots,x^{n})$, $\alpha=1,\dots,s$.

More precisely, letting   $\bar x^k=x^k(p)$, then in a neighborhood of $p$ the integral manifold $S$ is uniquely determined 
by the following  initial data:
\beq\label{ck-ic}
	\begin{array}{lclc}
 	f^{\alpha}&:=  & F^{\alpha}(\bar x^1,\bar x^2,\dots,\bar x^n)& \mbox{ for } {0}<\alpha\leq c_{0}, \\
 	f^{\alpha}(x^1)&:=  & F^{\alpha}(x^{1},\bar x^2,\dots,\bar x^n)& \mbox{ for } c_{0}<\alpha\leq c_{1},  \\
	f^{\alpha}(x^1, x^{2})&:=   & F^{\alpha}(x^{1},x^{2},\bar x^3,\dots, \bar x^n)& \mbox{ for } c_{1}<\alpha\leq c_{2},   \\
  	&  \vdots &   &\\
  	f^{\alpha}(x^1, x^{2},\dots,x^{n})&:=   & F^{\alpha}(x^{1},x^{2},\dots,x^{n})& \mbox{ for } c_{n-1}<\alpha\leq c_{n}.  
	\end{array}
\eeq
Namely,  for each $k=0,\dots,n$ and $\alpha$ such that $c_{k-1}<\alpha\leq c_k $,  let $f^\alpha(x^1,\dots,x^k)$ 
be an arbitrary analytic  function of $k$ variables, such that  $|f^\alpha(x^1,\dots,x^k)-y^\alpha(p)|$ is sufficiently 
small  in a neighborhood  of $(\bar x^1,\dots,\bar x^k)$
(where we define  $c_{-1}:= 0$ and by a function of zero variables we mean a constant). 
Then there exists a unique analytic integral manifold $y^{\alpha}=F^{\alpha}(x^{1},\dots,x^{n})$, $\alpha=1,\dots,s$ 
that satisfies the initial conditions  \eq{ck-ic}.
\end{theorem}

We now return to the $\lam$-system. To structure the presentation we first relabel indices 
(if necessary) and make the assumption 
\[\mbox{(A)}\quad\qquad\mbox{the unique algebraic relation \eq{alpha-lambda} 
does not involve $\lam^1$, i.e.\ $\alpha_1=0$.}\]  
\begin{theorem}\label{case_IIb}
	Assume $n=3$ and that the $\lambda$-system contains a single algebraic constraint 
	\eq{alpha-lambda} ($\rank(N)=1$). 
	Consider the case IIb where exactly two $\lambda^i$ appear with non-vanishing coefficients in 
	\eq{alpha-lambda}, and assume without loss of generality a labeling such that the 
	assumption (A) holds.
	
	If $\Gamma^3_{31}\neq \Gamma^2_{21}$ then the only solutions to the 
	$\lambda$-system are trivial solutions $\lam^1=\lam^2=\lam^3=constant$, while if 
	$\Gamma^3_{31}= \Gamma^2_{21}$ then the general solution depends on one arbitrary function of 
	one variable and one arbitrary constant. Finally, both situations can occur.
\end{theorem}
\begin{proof}
The algebraic relation \eq{alpha-lambda} is equivalent to \eq{c1}-\eq{c3}, and
from (A) it follows that 
\beq\label{alpha1=0} 
	c_{32}^{1}=\Gamma_{32}^1-\Gamma_{23}^1=0\,,\qquad \Gamma_{31}^2=0\,, \qquad\mbox{and} 
	\qquad\Gamma_{21}^3=0\,.
\eeq
Thus the algebraic relations \eq{c1}-\eq{c3} reduce to:
\begin{eqnarray*}
\Gamma_{32}^1(\lam^2-\lam^3)&=&0\\
\Gamma_{13}^2(\lam^2-\lam^3)&=& 0\\
\Gamma_{12}^3(\lam^2-\lam^3)&=& 0\,. 
\end{eqnarray*}
By assumption this system has rank 1 such that the algebraic equations in the $\lambda$-system
are satisfied if and only if $\lam^3= \lam^2$.
Using this relation to eliminate $\lam^3$ in the PDE system \eq{pde1}-\eq{pde6} we obtain:
\begin{eqnarray}
	r_1 (\lam ^2)&=&\Gamma^2_{21}(\lam^1-\lam^2) \label{pde11}\\ 
	r_1 (\lam^2)&=&\Gamma^3_{31}(\lam^1-\lam^2)  \label{pde12}\\
	r_2 (\lam^1)&=&\Gamma^1_{12}(\lam^2-\lam^1)  \label{pde13}\\
	r_2 (\lam^2)&=&0  \label{pde14}\\
	r_3 (\lam ^1)&=&\Gamma^1_{13}(\lam^2-\lam^1)  \label{pde15}\\
	r_3 (\lam ^2)&=&0\,.\label{pde16}
\end{eqnarray}
It follows immediately from the first two equations that if 
$\Gamma^3_{31}\neq \Gamma^2_{21}$, 
then $\lam^1=\lam^2=\lam^3$, and thus, by Proposition \ref{equal_lambdas}, 
the only solutions are trivial: $\lam^1=\lam^2=\lam^3\equiv constant$.
On the other hand, if $\lam^3=\lam^2$, such that the algebraic equations in 
the $\lambda$-system are satisfied, and in addition
\beq\label{N1-2-nontrivial}
	\Gamma^3_{31}= \Gamma^2_{21}\,,
\eeq
then the $\lambda$-system reduces to the sub-system \eq{pde12}-\eq{pde16} 
of five PDEs on $\RR^3$ for two unknowns $\lam^1$ and $\lam^2$. 
As before, $\lam^1= \lam^2\,\, (=\lam^3=constant)$ provides a trivial solution
to the $\lambda$-system by Proposition \ref{equal_lambdas}. Example~\ref{nonrich_IIb_trivial}
shows that such $\lam$-systems exist.

To analyze non-trivial solutions we now assume, in addition to \eq{alpha1=0} 
and \eq{N1-2-nontrivial}, that $\lam^1\neq \lam^2=\lam^3$. We note that the 
sub-system \eq{pde12}-\eq{pde16} is {\em not} of Frobenius type, and instead we 
apply the Cartan-K\"ahler theorem. This requires some preliminary calculations. 

We start by verifying the 
integrability conditions corresponding to equality of all mixed 2nd derivatives that 
{\em can} be computed from \eq{pde12}-\eq{pde16}:
$[r_3,r_2](\lam^1) = \sum c_{32}^k r_k(\lam^1)$, 
$[r_2,r_1](\lam^2) = c_{21}^k r_k(\lam^2)$, 
$[r_3,r_1](\lam^2) = \sum c_{31}^k r_k(\lam^2)$, and 
$[r_3,r_2](\lam^2) = \sum c_{32}^k r_k(\lam^2)$, where the summations are over 
$k=1,\, 2,\, 3$.
Taking into account  the assumptions  $\lam^1\neq \lam^2$ and $c_{23}^1=0$, 
these reduce to:  
\begin{eqnarray}
	r_3(\Gamma^1_{12})-r_2(\Gamma^1_{13})
	&=&c^2_{32}\Gamma_{12}^1+c^3_{32}\Gamma_{13}^1 \label{r23lambda1}\\ 
	r_2(\Gamma^2_{21})&=&\Gamma_{21}^2\Gamma^1_{21}  \label{r21lambda2}\\
	r_3(\Gamma^2_{21})&=&\Gamma^2_{21}\Gamma^1_{31}\,.\label{r31lambda2}
\end{eqnarray}
\begin{lemma}\label{0_torsion}
	Due to symmetry \eq{T0} and flatness \eq{R0}, and the
	assumptions \eq{alpha1=0} and \eq{N1-2-nontrivial}, the compatibility conditions 
	\eq{r23lambda1}, \eq{r21lambda2}, \eq{r31lambda2} are all satisfied as identities.
\end{lemma}
\noindent
{\em Proof of Lemma \ref{0_torsion}.} We only include the most involved case 
of \eq{r23lambda1}. Equation \eq{R0} with $m=3$, $k=1$, $i=2$, and $j=1$,
after simplifications due to \eq{T0}, reads:
\[r_3(\Gamma^1_{12})-r_1(\Gamma^1_{32})= 
\big(\Gamma^1_{11}\Gamma^1_{32}-\Gamma^1_{13}\Gamma^1_{12}\big)
+\big(\Gamma^1_{12}\Gamma^2_{32}-\Gamma^1_{32}\Gamma^2_{12}
-c^2_{13}\Gamma^1_{22} \big)
+\big(\Gamma^1_{13}\Gamma^3_{32}-\Gamma^1_{33}\Gamma^3_{12}
-c^3_{13}\Gamma^1_{32}\big)\, . \]
Similarly, equation \eq{R0} with $m=2$, $k=1$, $i=3$, $j=1$ gives:
\[r_2(\Gamma^1_{13})-r_1(\Gamma^1_{23})=\big(\Gamma^1_{11}\Gamma^1_{23}
-\Gamma^1_{12}\Gamma^1_{13}\big)
+\big(\Gamma^1_{12}\Gamma^2_{23}-\Gamma^1_{22}\Gamma^2_{13}
-c^2_{12}\Gamma^1_{23}\big)
+\big(\Gamma^1_{13}\Gamma^3_{23}-\Gamma^1_{23}\Gamma^3_{13}
-c^3_{12}\Gamma^1_{33}\big)\, .\]
Subtracting these and keeping in mind \eq{T0} and assumptions 
\eq{alpha1=0} and \eq{N1-2-nontrivial}, we get
\beann
	r_3(\Gamma^1_{12})-r_2(\Gamma^1_{13})&=&\left(\Gamma^1_{12}\Gamma^2_{32}
	-\Gamma^1_{23}\Gamma^2_{12}-c^2_{13}\Gamma^1_{22}+\Gamma^1_{13}\Gamma^3_{32}
	-\Gamma^1_{33}\Gamma^3_{12}-c^3_{13}\Gamma^1_{23}\right)\\
	&&-\,\,\left( \Gamma^1_{12}\Gamma^2_{23}-\Gamma^1_{22}\Gamma^2_{13}
	-c^2_{12}\Gamma^1_{23}+\Gamma^1_{13}\Gamma^3_{23}-\Gamma^1_{23}\Gamma^3_{13}
	-c^3_{12}\Gamma^1_{33}\right)\,\nn\\
	&=&\Gamma^1_{12}c^2_{32}-\Gamma^1_{23}\Gamma^2_{21}+
	\Gamma^1_{22}\Gamma^2_{31}+\Gamma^1_{13}c^3_{32}
	-\Gamma^1_{33}\Gamma^3_{21}+\Gamma^1_{23}\Gamma^3_{31}\\
	&=&\Gamma^1_{12}c^2_{32}+\Gamma^1_{13}c^3_{32}\,.
\eeann
which is \eq{r23lambda1}. Similar calculations verify \eq{r21lambda2} and \eq{r31lambda2}.
\qed

\medskip

Returning to the proof of Proposition \ref{case_IIb} we proceed to analyze the EDS 
$\mathcal I$ associated to the PDE system \eq{pde12}-\eq{pde16}, which is  
 differentially generated by the 1-forms
\begin{eqnarray*}
	\theta^1&=&d\lam^1-s\, \ell^1
	-\Gamma^1_{12}(\lam^2-\lam^1)\,\ell^2-\Gamma^1_{13}(\lam^2-\lam^1)\,\ell^3\\
	\theta^2&=&d\lam^2-\Gamma^2_{21}(\lam^1-\lam^2)\,\ell^1\,,
\end{eqnarray*}
on $M:=\RR^6$ with coordinates $u^1$, $u^2$, $u^3$, $\lam^1$, $\lam^2$, and $s$, where 
$s$ represents $r_1(\lam^1)$. Lemma \ref{0_torsion} amounts to the fact that the essential 
torsion of $\mathcal I$ vanishes identically. 
A direct computation (making use of Lemma \ref{0_torsion}, $c^1_{32}=0$ and 
$\Gamma^3_{31}=\Gamma^2_{21}$) shows that 
\beann
	d\theta^1&\equiv&\pi\wedge \ell^1 \quad \mbox{ mod } \{\theta^1,\theta^2\} \\
 	d\theta^2&\equiv&0 \quad  \mbox{ mod } \{\theta^1,\theta^2\}\,,
\eeann
where 
\begin{eqnarray} 
	\pi &=& -ds+\big[-s(c^1_{12}+\Gamma^1_{12})+(\lam^2-\lam^1)(r_1(\Gamma^1_{12})
	-\Gamma^1_{12}\Gamma^2_{12}-\Gamma^1_{13}c^3_{12})\big]\,\ell^2\nn \\
	&&\qquad+\,\big[-s(c^1_{13}+\Gamma^1_{13})+(\lam^2-\lam^1)(r_1(\Gamma^1_{13})
	-\Gamma^1_{13}\Gamma^3_{13}-\Gamma^1_{12}c^2_{13})\big]\,\ell^3\label{pi}\\
	&=:&-ds+A_2\ell^2+A_3\ell^3\,.\nn
\end{eqnarray}
Thus, the 1-forms $\ell^1,\ell^2,\ell^3, d\lam^1, d\lam^2, \pi$ provide a coframe on $\RR^6$, 
and the  EDS $\mathcal I$ is {\em algebraically} generated by $\theta^1, \theta^2$, and 
$\pi\wedge \ell^1$. 

The remaining parts of the proof consist in: (1) describing the variety of 3-dimensional 
integral elements of $\mathcal I$ that satisfy an independence condition, (2) choosing 
a flag of integral elements, and computing corresponding polar spaces, and (3) applying the 
Cartan test and the Cartan-K\"ahler theorem to determine the set of solutions.
   
(1) {\em The variety of integral elements:} let $E_3\in G_3(TM|_\point) $ be an integral element of $\mathcal I$ at the 
arbitrary point $\point\in M\cong\RR^6$, where  $G_3(TM|_\point)$ denotes 
the Grassmannian manifold of 3-dimensional subspaces of the tangent space $T_\point M$. 
The 3-plane $E_3$ is required to satisfy the independence condition
\beq
	\ell^1\wedge \ell^2\wedge \ell^3|_{E_3}\neq 0\,.\label{ind-ell3}
\eeq 
We let $\{e_1,\, e_2,\, e_3\}$ be a basis of $E_3$ such that 
\beq\label{ell-e3} 
	\ell^i(e_j)=\delta_{ij}\qquad i,j=1,\, 2,\, 3\,.
\eeq
The following are then necessary and sufficient conditions for $E_3$ to be an integral element of $\mathcal I$:
\[ \theta^1(e_i)=0,\qquad \theta^2(e_i)=0 \qquad\mbox{and}\qquad \pi\wedge \ell^1(e_i,e_j)=0\,, 
\qquad\mbox{ where }\, i,j=1,2,3\, . \]
Due to conditions \eq{ell-e3} these are equivalent to 
\beq\label{int-n3-rank1-2}
	\theta^1(e_i)=0\,, \qquad\theta^2(e_i)=0\,, \qquad \mbox{for}\quad i,j=1,2,3\quad \mbox{and} \qquad \pi(e_2)=0\,,
	\qquad \pi(e_3)=0\, .
\eeq
Thus the variety $\mathcal V _3(\mathcal I_\point) \subset G_3(TM|_\point) $ of 
3-dimensional integral elements is defined by 8 independent linear equations \eq
{int-n3-rank1-2}. We conclude that 
\beq
\nonumber	\codim \mathcal V _3(\mathcal I|_\point)=8,  
	\mbox{ and so } \dim \mathcal V _3(\mathcal I|_\point)
	=\dim G_3(TM|_\point)-\codim \mathcal V _3(\mathcal I|_\point)=9-8=1\,.
\eeq    
Explicitly  we can write 
\beq
	\nonumber e_i=r_i+p^1_i\dbyd{\lam^1}+p^2_i\dbyd{\lam^2}+a_i\dbyd{s}\,,\qquad \mbox{for}\quad i=1,2,3\, ,
\eeq 
where $p_i^1, p_i^2, a_i,\, i=1,\, 2,\, 3,$ are coordinate functions on an open subset of 
$G_3(TM|_\point)$ on which the independence condition \eq{ind-ell3} is satisfied. 
The conditions \eq{int-n3-rank1-2} imply that $\mathcal V_3(\mathcal I)$ is defined by the equations 
\beq
	\label{vi3} 
	p^1_1=s\,,\qquad
	p^1_2=\Gamma^1_{12}(\lam^2-\lam^1)\,,\qquad
	p^1_3=\Gamma^1_{13}(\lam^2-\lam^1)\,,\qquad
	p^2_1=\Gamma^2_{21}(\lam^1-\lam^2)\,,
\eeq
and
\beq
	p^2_2=0\,,\qquad
	p^2_3=0\,,\qquad
	a_2=A_2\,,\qquad
	a_3=A_3\,,
\eeq
where $A_2$ and $A_3$ are defined in \eq{pi}. The variety
$\mathcal V_3(\mathcal I)$ is therefore parametrized by a unique coordinate function 
$a_1$. Equations \eq{vi3} are of constant rank at every point $\point\in M$, and therefore all 
integral elements are ordinary at every point of $M$.

(2){\em Flag of integral elements:}  
By setting the arbitrary parameter $a_1=0$ we specify a \emph{particular} integral plane 
$E_{3}=\spa\{e_1,e_{2}, e_3\}$, where
\bes
	\nonumber  e_1&=&r_1+s\dbyd {\lam^1}+\Gamma^2_{21}(\lam^1-\lam^2)\dbyd {\lam^2}\,,\\
	\nonumber  e_2&=&r_2+\Gamma^1_{12}(\lam^2-\lam^1)\dbyd {\lam^1}+A_2\dbyd {s}\,,\\
	\nonumber  e_3&=&r_3+\Gamma^1_{13}(\lam^2-\lam^1)\dbyd {\lam^1}+A_3\dbyd {s}\,.
\ees
We define a flag of integral elements:
\beq
	\label{flag-n3-rank1-2} E_0=\{0\}\subset E_1=\spa\{e_1\}\subset 
	E_2=\spa\{e_1,e_2\}\subset E_3=\spa\{e_1,e_2,e_3\}\,, 
\eeq
A calculation, using that $\pi(e_1)=0$, $\ell^1(e_1)=1$, and  $\ell^1(e_2)=0$, shows that the corresponding polar spaces are
\bes
	\nonumber  H(E_0)&=&\{\theta^1, \theta^2\}^\perp,\\
 	\nonumber  H(E_1)&=&\{v\in H(E_0)|\pi\wedge 
	\ell^1(e_1,v)=0\}=\{\theta^1, \theta^2,\pi\}^\perp\\
	\nonumber  H(E_2)&=&\{v\in H(E_0)|\pi\wedge \ell^1(e_1,v)=0
	\quad \mbox{and} \quad \pi\wedge \ell^1(e_2,v)=0 \}=\{\theta^1, \theta^2, \pi\}^\perp \,.
\ees
Therefore
\beq 
	\nonumber c_0:=\codim H(E_0)=2\,,\qquad c_1:=\codim H(E_1)=3\,,\qquad c_2:=\codim H(E_2)=3\,.
\eeq
(3) {\em Applying Cartan test and the Cartan-K\"ahler theorem:} 
 Since $c_0+c_1+c_2=\codim\mathcal V _3(\mathcal I)=8$ the equality in the Cartan's Test for involutivity holds. Therefore
each $E_k$ in  the flag \eq {flag-n3-rank1-2}  is K\"ahler regular (see Theorem 7.4.1 \cite{il}).  

According to the Cartan-K\"ahler Theorem~\ref{Cartan-Kahler} the general solution depends on 2 constants that prescribe the values 
of $\lam^1$ and $\lam^2=\lam^3$ at an initial point $\bar u$ and one arbitrary function of one variable that prescribes the directional derivative  $s=r_1(\lam^1)$ along a curve. This arbitrary function absorbs the arbitrary constant that     prescribe the value
of $\lam^1$ at $\bar u$ and thus effectively  the general solution depends on one arbitrary function of one variable and one arbitrary constant.
The $\lam$-system in Examples~\ref{nonrich_IIb_non-trivial}
 has solutions of this type. This concludes the proof of Proposition \ref{case_IIb}.
\end{proof}

%
%
\section{Rich systems}\label{rich}
In this section we consider Problem \ref{prob1} for rich systems where we are given a 
coordinate system $u$ and a frame $\{R_1(u),\dots,R_n(u)\}$ satisfying Definition 
\ref{rich_sys}. This corresponds to searching for systems \eq{claw} whose eigencurves 
are the coordinate curves of some system of coordinates $w^1(u),\dots,w^n(u)$
(the Riemann invariants) on $\RR^n$. This leads to a slight reformulation of 
Problem \ref{prob1}, and a theorem of Darboux provides a complete solution 
for any dimension $n$. 

The reason why general rich systems can be completely analyzed in this way is the
fact that all the algebraic constraints in this case either are absent, or always 
impose equality of pairs of eigenvalues. This is in contrast to general (non-rich)
systems for which the algebraic conditions may be more complicated; see 
Section \ref{rank1_n=3}.

\subsection{The $\lambda$-system in Riemann invariants} 
We recall that richness may be formulated in several equivalent ways. (For the 
setting where a system \eq{claw} is given, see Chapter 12 in \cite{serre2} 
or Sections 7.3-7.4 in \cite{daf}.)  
In particular, it follows from Definition \ref{rich_sys} and the Frobenius theorem 
that the frame $\{R_1(u),\dots,R_n(u)\}$ is 
rich if and only if there is a change of coordinates 
\[u\mapsto \rho(u)=(w^1(u),\dots,w^n(u))\qquad\mbox{with}\qquad \nabla w^i(u)\cdot R_j(u)  \left\{\begin{array}{ll}
	=0 & \mbox{if $i\neq j$}\,,\\
	\neq 0 & \mbox{if $i=j$}\,.
\end{array}\right.\] 
The $w$-coordinates are  referred to as associated {\em Riemann invariants}. 
These are not unique and we assume that we have fixed one choice of the map $w=\rho(u)$.
In this case we may scale the given vector fields $R_i(u)$ according to the normalization
\beq\label{norm}
	\nabla w^i(u)\cdot R_j(u)\equiv \delta^i_j\,,
\eeq
which we assume throughout this section. A calculation (see \cite{daf} Section 7.3) then shows that 
these scalings render $\{R_1(u),\dots,R_n(u)\}$ a {\em commutative} frame: all structure 
coefficients $c^i_{jk}$ in \eq{lie-bracket} vanish. 
Furthermore, according to the normalization \eq{norm}, the corresponding left eigenvectors 
$L^i(u)$ in \eq{right_left}${}_2$ 
are given by $L^i(u)=\nabla w^i(u)$, such that the matrix $L(u)$ 
is a Jacobian with respect to $u$: $Ldu=dw$.  

In $w$-coordinates the  $\lambda$-system  \eq{sev1} - \eq{sev2} becomes
\begin{eqnarray}
	\del_i \ka^j &=& Z_{ji}^j (\ka^i-\ka^j)
	\qquad \mbox{for $1\leq i\neq j\leq n$,}\qquad \quad \big(\del_i=\textstyle\pd{}{w^i}\big)\label{sev1rich}\\
	Z_{ij}^k (\ka^j-\ka^i)&=&0 
	\qquad\qquad\qquad\,\,\,\, \mbox{for $1\leq k\neq i\,<j \,\neq k\leq n$,}\label{sev2rich}
\end{eqnarray}
where
\beq\label{pull_backs}
	 \ka^i(w):=\lam^i\circ\rho^{-1}(w)
	\quad\mbox{and} \quad Z_{ij}^k(w):=\Gamma_{ij}^k\circ\rho^{-1}(w)\,.
\eeq
Symmetry \eq{T0} and flatness \eq{R0} of the connection $\nabla$  imply  the following properties of 
the Christoffel symbols $Z_{ij}^k(w)$:
\beq\label{T0Z} 
	Z_{ij}^k= Z_{ji}^k\qquad\mbox{ (symmetry), }
\eeq
\beq\label{R0Z} 
	\del_m\big(Z^j_{ik}\big)-\del_k\big(Z^j_{im}\big) 
	=\sum_{t=1}^n\big(Z^j_{tk}Z^t_{im}-Z^j_{tm} Z^t_{ik}\big)\,,\qquad\forall\, i,\, j,\, k,\, m\,\mbox{ (flatness)}.
\eeq
Problem \ref{prob1} in the rich case thus takes the following form:
\begin{problem}\label{rich_prob}(Rich frame) 
	With the same notation as in Section \ref{statement}, assume that the given frame 
	$\{R_1(u),\dots,R_n(u)\}$ is rich, and let $w=(w^1,\dots,w^n)$ be associated 
	Riemann invariants. Assume the normalization \eq{norm}, and define the connection 
	coefficients $Z_{ij}^k$ by \eq{pull_backs}. Then, determine the set of solutions 
	$\ka^1(w),\dots,\ka^n(w)$ of the $\lambda$-system in Riemann invariants
	\eq{sev1rich}-\eq{sev2rich}.
\end{problem}
Below we give a complete answer to this problem by applying a theorem of 
Darboux. For completeness we include a precise statement of this result.

\begin{theorem}\label{dar3} (Darboux \cite{dar}) Given a system of first order PDEs 
with dependent variables $v=(v^1,\dots,v^m)$ and independent variables $x=(x^1,\dots,x^n)$. 
Assume that each equation is of the form
\beq\label{sys_ii}
	\frac{\del v^i}{\del x^h} = f_{ih}(v,x)\,,
\eeq
where the given maps $v\mapsto f_{ih}(v,x)$ are $C^1$-smooth, uniformly 
in $x$ in a neighborhood of a given point $x_0=(x^1_0,\dots,x^n_0)$. 
Assume that the system prescribes compatible second order mixed derivatives 
in the following sense:
\begin{itemize}
	\item[(C)] Whenever $v^i$ is such that two different first derivatives 
	$\frac{\del v^i}{\del x^h}$ and $\frac{\del v^i}{\del x^{k}}$ are given by the system 
	\eq{sys_ii}, the equation
	\[\frac{\del}{\del x^{k}}\big(f_{ih}(v,x)\big) = \frac{\del}{\del x^{h}}\big(f_{ik}(v,x)\big)\]
	(after taking the total derivative on each side),
	contains only first order derivatives which are prescribed by \eq{sys_ii}, and
	substitution of \eq{sys_ii} for the first derivatives that appear gives
	an identity in $x$ and $v$.
\end{itemize} 
Suppose a dependent variable $v^i$ appears differentiated in \eq{sys_ii} 
with respect to $x^{h_1},\dots,x^{h_p}$. Then, letting $\tilde x$ denote the remaining 
independent variables, we prescribe a real-valued $C^1$-smooth function $\varphi^i(\tilde x)$ and require that
\beq\label{DDdata}
	v^i(x^1,\dots,x^n)\Big|_{x^{h_1}=x^{h_1}_0,\dots,\, x^{h_p}=x^{h_p}_0}
	=\varphi^i(\tilde x)\,.
\eeq
We make such an assignment of data for each dependent variable that appears 
differentiated in \eq{sys_ii}. Then, under the compatibility condition (C), the problem 
\eq{sys_ii} - \eq{DDdata} has a unique local $C^2$-solution.
\end{theorem}
\begin{remark}
	The initial data \eq{DDdata} provide the values of each 
	dependent variable in 
	the directions where its derivatives are not prescribed by the system \eq{sys_ii}. 
	In particular, if all partial derivatives of a variable $v^j$ are prescribed by \eq{sys_ii}, 
	then $\varphi^j(\tilde x)$ is taken to be a constant.
	
	Our formulation of Theorem \ref{dar3} covers several situations of both over-
	and under-determined systems. In \cite{dar} the different cases are stated in three
	separate theorems, one of which is the PDE version of the Frobenius theorem (all 
	first partials of all dependent variables are prescribed). The theorem is proved by 
	a suitable Picard iteration and requires only $C^1$-regularity.
\end{remark}

%
%
\subsection {Rich systems with no algebraic constraints}\label{rich-rank0}
We first consider the situation when $Z^k_{ij}=0$ for {\em all} triples of distinct $i,\, j,\, k$, 
such that the algebraic constraints \eq{sev2rich} are trivially satisfied, i.e., $\rank(N)=0$.
We verify that the compatibility conditions (C) in Theorem \ref{dar3} are satisfied 
in this case, and that the general solution of the system \eq{sev1rich} depends on $n$ arbitrary 
functions of one variable. In particular, there are many strictly hyperbolic solutions of 
Problem \ref{rich_prob} in this case. 

We also show that if the given vector fields $R_i$ 
(in addition to forming a rich frame) are orthogonal, 
then the corresponding $\lambda$-system necessarily belongs to this unconstrained case. 
Two concrete examples of rich and algebraically unconstrained systems are the class of 
$2\times 2$-systems, and the class of  $n\times n$-systems
with constant eigenfields, see Examples \ref{n2}, \ref{rich_orth}, and \ref{constant_eigenfields}.
\begin{theorem}\label{rich_unconstr} 
	Given a $C^1$-smooth, rich frame $\{R_1(u),\dots,R_n(u)\}$ in a neighborhood 
	of $\br{w}\in\RR^n$. Let $\rho(u)=(w^1(u),\dots,w^n(u))$ be associated Riemann 
	invariants and assume the normalization \eq{norm}.
	Let the connection coefficients $Z_{ij}^k$ be defined by \eq{pull_backs} and 
	assume that $Z^k_{ij}=0$ whenever $i\neq j\neq k\neq i$. 
	Then, for given functions $\varphi_i$, $i=1,\dots,n$, of one variable, there is a unique local 
	solution $\ka^1(w),\dots,\ka^n(w)$ to the $\lambda$-system 
	\eq{sev1rich} with
	\[\ka^i(\br{w}^1,\dots,\br{w}^{i-1}, w^i, \br{w}^{i+1}, \dots \br{w}^n)=\varphi_i(w^i).\]
\end{theorem}
\pf  The $\lambda$-system \eq{sev1rich} consists of $n(n-1)$ PDEs that prescribe, for each $j=1,\dots,n$,
all first partials of $\kappa^j$ except $\del_j\kappa^j$. 
The conclusion of the theorem thus follows from Darboux's Theorem \ref{dar3} provided the system \eq{sev1rich} 
satisfies the compatibility conditions (C). That is, for each $ j=1,\dots, n$ the equalities 
\[\del_k\del_m\ka^j=\del_m\del_k \ka^j\]
should hold as identities for all $k \neq j$, $m\neq j$ and $m\neq k$, when the first derivatives of $\ka$'s 
are substituted for from \eq{sev1rich}. This leads to the compatibility conditions
\begin{eqnarray}\label{int-k0}  
	\left(\del_{m}Z^{j}_{jk} -\del_{k}Z^{j}_{jm}\right)\ka^j
	&+&\left(Z^{j}_{jm}Z^{m}_{mk}+Z^{j}_{jk}Z_{km}^{k}-Z^{j}_{jm}Z^{j}_{jk}-\del_{m}Z^{j}_{jk}\right)\ka^{k}\\
	\nonumber&-&\left(Z^{j}_{jm}Z^{m}_{mk}+Z^{j}_{jk}Z_{km}^{k} - Z^{j}_{jk}Z^{j}_{jm}
	-\del_{k}Z^{j}_{jm}\right)\ka^m\equiv 0\,,
\end{eqnarray}
for all distinct $j,\, m,\, k$.
We verify these conditions by showing that the coefficients of $\ka^j$, $\ka^m$ and $\ka^k$ 
vanish identically due to \eq{T0Z} and the flatness condition \eq{R0Z}. We first substitute  $i=j$ 
in \eq{R0Z} to obtain
\[\del_m Z^j_{jk}-\del_kZ^j_{jm}=\sum_{t=1}^n \big(Z^j_{tk}Z^t_{jm}-Z^j_{tm}Z^t_{jk}\big)\,,\]
which vanishes since $k,\, m,\, j$ are distinct and since, by assumption, $Z^{i_1}_{i_2i_3}\equiv 0$ for 
distinct $i_1,\, i_2,\, i_3$. This shows that the coefficient of $\ka^j$ in \eq{int-k0} is zero. 
The arguments for $\ka^k$ and $\ka^m$ are similar and we only consider the coefficient 
of $\ka^k$. Interchanging $k$ and $i$ in \eq{R0Z} and setting $i=j$ yields
\[\del_m Z^j_{jk}-\del_jZ^j_{mk}=\sum_{t=1}^n \big(Z^j_{jt}Z^t_{mk}-Z^j_{mt}Z^t_{jk}\big)\,.\]
Again, using \eq{T0Z}, that $Z^{i_1}_{i_2i_3}\equiv 0$ for distinct $i_1,\, i_2,\, i_3$, and that 
$k,\, m,\, j$ are distinct, we obtain
\[\del_m Z^j_{jk}=Z^j_{jm}Z^m_{mk}+Z^j_{jk}Z^k_{mk}-Z^j_{mj}Z^j_{jk}\,.\] \foorp

\subsubsection{\bf{Rich, orthogonal frames}}\label{rich-orthogonal}
Assume now that the given rich frame $\{R_1(u),\dots,R_n(u)\}$ has the additional 
property that it is orthogonal: $R_i\cdot R_j=0$ if $i\neq j$. That is, we search for systems
\eq{claw} whose eigencurves are the coordinate curves of an {\em orthogonal coordinate system}
$(w^1(u),\dots,w^n(u))$ on $\RR^n$. In this case we show that the connection components 
$Z^k_{ij}$ necessarily vanish whenever $i,\, j,\, k$ are distinct, such that Theorem 
\ref{rich_unconstr} applies.

We define a matrix $S(w)$, whose components $S^i_j(w):=R^i_j\circ\rho^{-1}(w)$ are the pull-backs 
of components of the matrix $R$ under $\rho^{-1}$. The connections one-forms (see \eq{conn-forms} and \eq {mu}) are then given by 
:
\beq\label{tmu}\mu_i^j:=\sum_{k=1}^nZ^j_{ki}\, dw^k=\big(S^{-1}dS\big)^j_i\,.\eeq
By orthogonality we have that 
\beq\label{orth}
	S^TS=\diag[|S_1|^2\dots|S_n|^2] \qquad\mbox{such that}\qquad
	S^{-1}=\diag\left[|S_1|^{-2}\dots|S_n|^{-2}\right]S^T\,.
\eeq 
It follows that the connection matrix $\mu$ is given by
\[\mu=S^{-1}dS=\diag\left[|S_1|^{-2}\dots|S_n|^{-2}\right]S^T\,dS.\]
Differentiation of \eq{orth} gives
\[ d(S^T)S+S^TdS=\diag[d|S_1|^2\dots d|S_n|^2]\,,\]
such that
\[\mu^T=
\diag\left[\frac {d|S_1|^2}{|S_1|^2}\dots\frac{d |S_n|^2}{|S_n|^2}\right]
-\diag\left[|S_1|^2\dots|S_n|^2\right]\,\mu\, \diag\left[|S_1|^{-2}\dots|S_n|^{-2}\right]\,.\]
Component-wise we thus have
\beq\label{tmuij}
	\mu^i_j = -\frac{|S_j|^2} {|S_i|^2}\mu^j_i 
	\qquad \forall j\neq i \qquad \mbox{and}\qquad 
	\mu^i_i = \frac 1{2\,|S_i|^2}d|S_i|^2\,. 
\eeq
From  \eq{tmu} and \eq{tmuij} it now follows that 
\[Z^{i}_{jl}\, =\, -\frac{|S_j|^2}{|S_i|^2}Z^j_{il} 
\qquad \forall l,\, \forall j\neq i\,,\]
By symmetry in the lower indices \eq{T0Z} we get, for indices $i\neq j\neq l\neq i$, 
that 
\beq\label{Z1}
	Z_{jl}^{i}=-\frac{|S_j|^2}{|S_i|^2}Z^{j}_{il}=-\frac{|S_j|^2}{|S_i|^2}Z^{j}_{li}
	=\frac{|S_l|^2}{|S_i|^2}Z^l_{ji}\,,
\eeq
while, at the same time,
\beq\label{Z2}
	Z^{i}_{jl} = Z^{i}_{lj} =  -\frac{|S_l|^2}{|S_i|^2}Z^l_{ij} = -\frac{|S_l|^2}{|S_i|^2}Z^l_{ji}\,. 
\eeq
It follows from \eq{Z1} and \eq{Z2} that $Z^{i}_{jl}=0$ whenever all three indices are distinct, such that Theorem 
\ref{rich_unconstr} applies. We conclude that for rich, orthogonal frames the solution of Problem \ref{rich_prob}
depends on $n$ functions of one variable.

\bigskip

\subsection{Rich systems with algebraic constraints}\label{rich-rank-non-zero}
We next consider the more involved situation of rich systems with non-trivial algebraic 
constraints \eq{sev2rich} such that $\rank(N)>0$: there exist distinct $i,\, j,\, k$ such that $Z^k_{ij}\neq 0$. 
The algebraic relations \eq{sev2rich} then impose the equality $\ka^i=\ka^j$,
such that there are multiplicity conditions on the eigenvalues. 
(In particular, there are no strictly hyperbolic systems in this ``rich \& $\rank(N)>0$" case.)
This makes the analysis of the $\lambda$-system more involved.

The analysis is further complicated by the fact that there may be additional 
equalities among eigenvalues that are imposed by the PDEs \eq{sev1rich}.
We claim that, after taking all of these relations into account, the $\lambda$-system 
\eq{sev1rich}-\eq{sev2rich} reduces to a first order system of PDEs to which 
the Darboux theorem (Theorem \ref{dar3}) can be applied. We have:

\begin{theorem} \label{rich_cnstrnd}
	Given a $C^1$-smooth, rich frame $\{R_1(u),\dots,R_n(u)\}$ in a neighborhood 
	of $\br{w}\in\RR^n$. Let $(w^1(u),\dots,w^n(u))$ be associated Riemann invariants
	and assume the normalization \eq{norm}.
	Let the connection coefficients $Z_{ij}^k$ be defined by \eq{pull_backs}, and 
	assume that there exists at least one triple of distinct $i,\, j,\, k$ with $Z^k_{ij} \neq 0$. 
	
	Then the $\lambda$-system \eq{sev1rich}-\eq{sev2rich} imposes multiplicity conditions 
	on the eigenvalues in the following sense. There are disjoint subsets 
	$A_1, \dots, A_{s_0}\subset \{1,\dots,n\}$ ($s_0\geq1$) of cardinality two or more, and such 
	that \eq{sev1rich}-\eq{sev2rich} impose the equality $\ka^i=\ka^j$ if and only if 
	$i,\, j\in A_\alpha$ for some $\alpha\in\{1,\dots,s_0\}$. 
	Let $l=\sum_{\alpha=1}^{s_0} |A_\alpha|\leq n$ and $s_1=n-l$. By relabeling indices we 
	may assume that $\{1,\dots,n\}\setminus\bigcup_{\alpha=1}^{s_0}A_\alpha=\{1,\dots,s_1\}$.

	Then, given $s_1$ functions of one variable $\varphi^1,\dots,\varphi^{s_1}$, and $s_0$ constants 
	$\hat\ha^1,\dots,\hat\ha^{s_0}$, there is a unique local solution $\kappa^1,\dots,\kappa^n$ of 
	the $\lambda$-system \eq{sev1rich}-\eq{sev2rich} such that
	\begin{eqnarray*}
 		\ka^i(\br{w}^1,\dots,\br{w}^{i-1}, w^i,\br{w}^{i+1}\dots \br{w}^n) &=& 
		\varphi_i(w^i) \,\,\qquad \mbox{if $i=1,\dots,s_1$}\,,\\
		\ka^i(\br{w}) &=& \hat\ha^\alpha \qquad\qquad \mbox{if $i\in A_\alpha$ for some $\alpha =1,\dots, s_0$.}
	\end{eqnarray*}
  	Moreover, $\del_i\ka^j(w)=0$ whenever $i,\, j\in A_\alpha$ for some $\alpha=1,\dots, s_0$.
	Here $s_1$ is the maximal number of simple eigenvalues, while $s_0$ 
	is the maximal number of non-simple eigenvalues, in a solution of the $\lambda$-system.
\end{theorem}
This result is a direct consequence of the following two lemmas and Theorem \ref{dar3}.
The first lemma details how the index sets $A_\alpha$ are defined, and how to write the  
$\lambda$-system \eq{sev1rich}-\eq{sev2rich} as a system of only PDEs.
\begin{lemma} \label{rich_cnstrnd1}
	With the same assumptions as in Proposition \ref{rich_cnstrnd}, the $\lambda$-system 
	\eq{sev1rich} - \eq{sev2rich} can be re-written as a first order PDE system without algebraic constraints. 
	More precisely, there are integers $s_0$, $s_1$, and disjoint subsets $A_1, \dots, A_{s_0}\subset \{1,\dots,n\}$ 
	as described in Proposition \ref{rich_cnstrnd} such that the following holds.
	The re-written system involves $s_1+s_0$ unknowns, denoted $\ka^1(w),\dots, \ka^{s_1}(w)$, 
	$\ha^1(w),\dots,\ha^\s0(w)$, and has the following form:
	\begin{eqnarray}
		\label{ki}
		\mbox{ for } j=1,\dots, s_1:&&\del_i{\ka^j}= \left\{\begin{array}{ll}
		Z_{ji}^j(\ka^i-\ka^j)&   \qquad \mbox{if $1 \leq i\neq j\leq s_1$} \\ 
		Z_{ji}^j(\ha^\alpha-\ka^j)&  \qquad  \mbox{if $i\in A_\alpha$ for some $\alpha =1,\dots, s_0$},
		\end{array}\right. \\ \nn\\
		\label{hi}
		\mbox{ for } \alpha=1,\dots, s_0: && \del_i{\ha^\alpha}=\left\{\begin{array}{ll}
		W_{i}^\alpha(\ka^i-\ha^\alpha)&  \qquad   \mbox{if $i=1,\dots,s_1$} \\ 
		W_{i}^\alpha(\ha^\beta-\ha^\alpha)& \qquad \mbox{if $i \in A_\beta$ and $1\leq \beta\neq\alpha\leq s_0$}\\
		0 &  \qquad \mbox{if $i\in A_\alpha$}\,.\end{array}\right.
	\end{eqnarray}
	Here the coefficients $W_i^\alpha$ are defined by $W_i^\alpha:=Z_{ji}^j$ if $j\in A_\alpha$ and $i\not\in A_\alpha$,
	and these are well-defined by the properties of the sets $A_\alpha$.
\end{lemma}
\pf
We construct the index sets $A_1,\dots, A_{s_0}$ by first considering the algebraic conditions \eq{sev2rich},
and then  taking into account the differential equations \eq{sev1rich}.

Each of the algebraic relations in \eq{sev2rich} is either trivial (``$0=0$")
or non-trivial, imposing equality of two eigenvalues. More precisely, 
if $1\leq i\neq j\leq n$ are such that  there exists 
$k\notin\{i,j\}$ with $Z^k_{ij}\neq 0$, then necessarily $\ka^i=\ka^j$.
By assumption there is at least one such non-trivial algebraic relation.
Grouping together the indices $i$ of the unknowns $\kappa^i$  
that must be identical according to these relations, we obtain a certain number $\ts0\geq1$ 
of disjoint index sets $\TA_1, \dots, \TA_{\ts0}\subset \{1,\dots,n\}$: two distinct $j_1,\, j_2$ 
belong to the same $\TA_\alpha$ if and only if one of the relations in \eq{sev2rich} imposes 
the equality $\ka^{j_1}=\ka^{j_2}$. Clearly, the cardinality of each $\TA_\alpha$ 
is at least two. For $\alpha=1,\dots, \ts0$ we introduce the unknowns $\tha^\alpha$ by setting
\[\tha^\alpha(w):= \ka^j(w) \qquad  \forall j\in \TA_\alpha\,,\]
which is well-defined by the definition of $\TA_\alpha$. Also let 
\beq\label{TA}
	\TA:=\bigcup_{\alpha=1}^{\ti s_0} \TA_\alpha \qquad\qquad \ti s_1:=n- |\TA|\,.
\eeq
If necessary we relabel the unknowns such that $\ka^1,\dots, \ka^{\ti s_1}$ 
denote the unknowns for which the algebraic relations \eq{sev2rich} impose no 
multiplicity constraint. 

We next turn to the PDE system \eq{sev1rich}. For each pair $(i,j)$, with $i\neq j$, either none, 
exactly one, or both of $i$ and $j$ belong to $\TA$. We list the possible cases together with the 
corresponding form of the PDEs in each case:
\beq\label{five_cases}
	\left\{\begin{array}{lll}
	\mathrm{(a)} & 1\leq i\neq j\leq \ti s_1: & \del_i \kappa ^j =Z_{ji}^j(\kappa ^i-\kappa ^j) \\\\
 	\mathrm{(b)} & \exists \alpha: \, i\in\TA_\alpha \mbox{ and } 1\leq j\leq \ti s_1 : 
	& \del_i \kappa ^j =Z_{ji}^j(\tha^\alpha-\kappa ^j) \\\\
	\mathrm{(c)} & \exists \beta: \, j\in \TA_\beta \mbox{ and } 1\leq i\leq \ti s_1 : 
	& \del_i \tha ^\beta =Z_{ji}^j(\kappa ^i-\tha^\beta) \\\\
	\mathrm{(d)} & \exists \gamma\neq \delta: \, i\in \TA_\gamma \mbox{ and } j\in \TA_\delta : & 
	\del_i \tha ^\delta =Z_{ji}^j(\tha ^\gamma-\tha ^\delta)\\\\
	\mathrm{(e)} &  \exists \epsilon: \, i,\, j\in \TA_\epsilon : & \del_i \tha ^\epsilon = 0\,.
	\end{array}\right.
\eeq
At this stage we have used all non-trivial algebraic relations imposed by the algebraic 
part \eq{sev2rich} of the $\lambda$-system. 
The issue now is that the cases (c) and (d) may impose further {\em algebraic} conditions
since their left-hand sides are independent of the index $j$:
\begin{itemize}
\item in case (c): for a given pair $(i,j)$ with 
$j\in \TA_\beta$ and $1\leq i\leq \ti s_1$, unless all the coefficients $Z_{ki}^k$
coincide as $k$ ranges over $\TA_\beta$,  we must impose that $\kappa ^i=\tha^\beta$.
If so we add $i$ to the index set $\TA_\beta$ and replace the $|\TA_\beta|$ PDEs
$\del_i \tha ^\beta =Z_{ji}^j(\kappa ^i-\tha^\beta)$ in (c) by the single PDE
$\del_i \tha ^\beta =0$. At the same time we substitute $\tha^\beta$ for $\kappa ^i$ in 
the remaining PDEs in (a)-(c) in which $\kappa^i$ appears.

On the other hand, if $Z_{ki}^k=Z_{ji}^j$ for all $k\in \TA_\beta$, then we define
$\ti W_i^\beta:=Z_{ji}^j$, and replace the same $|\TA_\beta|$ PDEs by the single PDE
$\del_i \tha ^\beta =\ti W_i^\beta(\kappa ^i-\tha^\beta)$.  
\item in case (d): for a given pair $(i,j)$ with $i\in \TA_\gamma$, $j\in \TA_\delta$,
and $\gamma\neq \delta$, unless the coefficients $Z_{ki}^k$ 
all coincide as $k$ ranges over $\TA_\delta$,  we must impose that $\tha ^\delta=\tha^\gamma$.
We then merge the index sets $\TA_\gamma$ and $\TA_\delta$ and replace the $|\TA_\delta|$ PDEs
$\del_i \tha ^\delta =Z_{ji}^j(\tha ^\gamma-\tha ^\delta)$ in (d) by the single PDE
$\del_i \tha ^\delta =0$. At the same time we substitute $\tha^\delta$ for $\tha ^\gamma$ in 
the remaining PDEs in (b)-(e) in which $\tha^\gamma$ appears.

On the other hand, if $Z_{ki}^k=Z_{ji}^j$ for all $k\in \TA_\delta$, then we define
$\ti W_i^\delta:=Z_{ji}^j$ and replace the same $|\TA_\delta|$ PDEs by the single PDE
$\del_i \tha ^\delta =\ti W_i^\delta(\tha ^\gamma-\tha ^\delta)$.
\end{itemize}
We continue this process of identifying unknowns $\kappa^i$ and $\tha^\alpha$ that must 
necessarily coincide (enlarging and merging index sets), until no further reduction 
is possible. Since each reduction consists in setting unknowns equal to each other, and since 
setting all unknowns equal gives the trivial solutions to the $\lambda$-system 
(Proposition \ref{equal_lambdas}), it follows that no contradiction will be obtained in this manner. 
Also, as there is a finite number of unknowns, the process must terminate after finitely many reductions.

At this point we have obtained a certain number $\s0\geq 1$ of disjoint index sets 
$A_1, \dots, A_{\s0}\subset \{1,\dots,n\}$ such that two distinct $j_1$, $j_2$ belong to $A_\alpha$ 
if and only if \eq{sev1rich} and \eq{sev2rich} impose the equality $\ka^{j_1}=\ka^{j_2}$ 
(through the reduction process described above). By assumption no further algebraic reduction is 
possible, and it follows that the sets $A_\alpha$ have the following properties:
\beq\label{Zassum1}  
	{Z^{j_1}_{j_1\,i}=Z^{j_2}_{j_2\,i}=:W^{\alpha}_{i}\qquad 
	\mbox{ whenever } j_1, j_2\in A_\alpha \mbox{ and } i\notin  A_\alpha} \,,
\eeq
and 
\beq\label{Zassum2} 
	{ Z^{k}_{ij}=0\qquad \mbox{ whenever }  
	i\in A_\alpha,\, j\notin A_\alpha \mbox{ and } k\notin\{i,j\} }\,.
\eeq
We also have that 
\beq\label{Zassum3} 
	{ Z^{k}_{ij}=0\qquad \mbox{ whenever }  
	1\leq i\neq j\leq s_1 \mbox{ and } k\notin\{i,j\} }\,.
\eeq
Each of the sets $\TA_\beta,\, \beta=1,\dots,\ts0$, 
is contained in one of the sets $A_\alpha,\, \alpha=1,\dots,\s0$, whence the latter have cardinality 
$|A_\alpha|\geq 2$. Setting
\[A:=\bigcup_{\alpha=1}^{s_0} A_\alpha\,,\qquad\mbox{and} \qquad s_1:=n-|A|\,,\]
it follows that $s_0$ is the maximal number of non-simple eigenvalues, and $s_1$ is the maximal 
number of simple eigenvalues, in a solution of the $\lambda$-system.
If necessary we relabel the indices such that $\kappa^1,\dots,\kappa^{s_1}$ denote the 
unknowns for which \eq{sev1rich}-\eq{sev2rich} impose no multiplicity constraints, and 
we introduce unknowns
\beq\label{h_kappa} 
	\ha^\alpha(w):= \ka^j(w) \qquad  \mbox{ for $j\in A_\alpha$, $\alpha=1,\dots, \s0\,.$}
\eeq
Similar to above, for each pair $(i,j)$ with $i\neq j$, either none, exactly one, or both of $i$ and $j$ belong
to $A$, and again there are five possible cases of the PDEs in \eq{sev1rich}. We thus obtain the same 
type of system as in \eq{five_cases} except with $\ti s_1$ replaced by $s_1$ and $\TA_\alpha$'s replaced 
by $A_\alpha$'s. Using \eq{Zassum1} we obtain the system \eq{ki}-\eq{hi}.
Finally, from the last equation in \eq{hi} it follows that $\ha^\alpha$ is independent of the coordinates 
$w^i$ for $i\in  A_\alpha$: $\del_i\ka^j(w)=0$ for all $i,\, j\in A_\alpha$.
\foorp
\begin{remark}
We note that if $s_0=1$ and $|A_1|=n$, i.e.\ $\ka^1=\dots=\ka^n=\ha^1$, then \eq{ki}-\eq{hi} 
degenerate into the trivial system $\pd{\ha^1}{w^i}=0,\,i=1,\dots,n$. In this case the only solutions
to the original problem are the trivial ones: $\ka^1(w)=\dots=\ka^n(w)\equiv \hat \ka\in\RR$. 
If $l=n$, then $s_1=0$ and the system \eq{ki} is empty. In this case all eigenvalues $\ka^i$ 
necessarily has  multiplicity greater than one.
\end{remark}

The PDE system \eq{ki} - \eq{hi} is in the form required by Theorem~\ref{dar3}, and we proceed 
to verify the compatibility conditions (C) in this case.

\begin{lemma} \label{int-hk-system}
	With the same assumptions as in Lemma \ref{rich_cnstrnd1}, the compatibility conditions (C) 
	in Theorem \ref{dar3} are satisfied for the rewritten $\lambda$-system \eq{ki}-\eq{hi}.
\end{lemma}
\pf \emph{Compatibility conditions for \eq{ki}:}
For $ j=1,\dots, s_1$, we need to verify that 
\beq\label{compat}
	\del_k\del_m\ka^j=\del_m\del_k \ka^j\qquad\mbox{for $k \neq j$, $m\neq j$ and $m\neq k$.}
\eeq
There are two main cases depending on whether $m$ and $k$ belong to a common 
$A_\alpha$ or not, and we first treat the latter. 

\noindent
{\em $\bullet$ Case 1.} $\nexists  \alpha$  such that $\{m,k\}\subset A_\alpha$. There are several 
sub-cases according to whether $m$ or $k$ belong to $A$. However, in each case 
the compatibility condition takes the form: 
\beq
	\begin{split}
	\left(\del_{m}Z^{j}_{jk} -\del_{k}Z^{j}_{jm}\right)\ka^j
	&+\left(Z^{j}_{jm}Z^{m}_{mk}+Z^{j}_{jk}Z_{km}^{k}-Z^{j}_{jm}Z^{j}_{jk}
	-\del_{m}Z^{j}_{jk}\right)\ka^{k}\label{int-k1}  \\ 
	&-\left(Z^{j}_{jm}Z^{m}_{mk}+Z^{j}_{jk}Z_{km}^{k} 
	- Z^{j}_{jk}Z^{j}_{jm}-\del_{k}Z^{j}_{jm}\right)\ka^m\equiv 0,
	\end{split}
\eeq
with the provision that 
\[\ka^m=\ha^\alpha\,,\qquad Z^m_{mk}=W^\alpha_k \qquad
\mbox{if $m \in A_\alpha$ for some $1\leq\alpha\leq \s0$,}\]
and
\[\ka^k=\ha^\beta\,,\qquad Z^k_{km}=W^\beta_m\qquad
\mbox{if $k \in A_\beta$ for some $1\leq\beta\leq \s0$.}\]
We verify that \eq{int-k1} is satisfied by showing that coefficients of $\ka^j$, $\ka^m$ 
and $\ka^k$ vanish identically due to the flatness condition \eq{R0Z}. 
We first substitute  $i=j$ in \eq{R0Z} to obtain
\[\del_m \big(Z^j_{jk}\big)-\del_k \big(Z^j_{jm}\big)
=\sum_{t=1}^n\big(Z^j_{kt}Z^t_{mj}-Z^j_{mt} Z^t_{kj}\big)\,.\]
By assumption \eq{Zassum2}, for distinct $k$, $m$, $j$ with $j\notin A$, we have 
\begin{eqnarray*}
	Z^t_{mj}&=&0\qquad\mbox{ unless }t=m\mbox{ or } t=j\\
	Z^t_{kj}&=&0\qquad\mbox{ unless }t=k\mbox{ or }t=j\,, 
\end{eqnarray*}	
such that
\beq\label{R0Zj=i} 
	\del_m \big(Z^j_{jk}\big)-\del_k \big(Z^j_{jm}\big)=Z^j_{mk}\big(Z^m_{mj}-Z^k_{kj}\big)\,.
\eeq
By assumption $\nexists \alpha$ such that $\{m,k\}\subset A_\alpha$, whence 
$Z^j_{mk}=0$ due to assumption \eq{Zassum2}. Thus
\beq\label{R0Zj=i1}
	\del_m \big(Z^j_{jk}\big) - \del_k\big(Z^j_{jm}\big) = 0\,,
\eeq
such that the coefficient of $\ka^j$ in \eq{int-k1} vanishes. 

Next consider the coefficient of  $\ka^k$ in \eq{int-k1}. Interchange $k$ and 
$i$ in \eq{R0Z} and set  $i:=j$ to get
\beq\label{flat}
	\del_m \big(Z^j_{jk}\big)-\del_j\big(Z^j_{mk}\big)=\sum_t\big(Z^j_{jt}Z^t_{mk}-Z^j_{mt}Z^t_{jk}\big)\,.
\eeq
Again, by assumption $\nexists \alpha$ such that $\{m,k\}\subset A_\alpha$, and \eq{Zassum2} gives
\beq\left\{\quad\begin{aligned}
	Z^j_{mk}&= 0\\
	Z^t_{mk}&= 0\qquad\mbox{ unless }t=m\mbox{ or } t=k\label{Zs}\\
	Z^t_{jk}&= 0\qquad\mbox{ unless }t=j\mbox{ or }t=k\,.
	\end{aligned}\right.
\eeq
Thus
\[\del_m Z^j_{jk}\,\,=\,\,
Z^j_{jm}Z^m_{mk} + Z^j_{jk}Z^k_{km} - Z^j_{jm}Z^j_{jk}-Z^j_{mk}Z^k_{kj}
\,\,=\,\,Z^j_{jm}Z^m_{mk} + Z^j_{jk}Z^k_{km} - Z^j_{jm}Z^j_{jk}\,,\]
such that the coefficient of $\ka^k$ in \eq{int-k1} vanishes.
The coefficient of $\ka^m$ is treated similarly.

\medskip

{\em $\bullet$  Case 2.} $\exists  \alpha$  such that $\{m,k\}\subset A_\alpha$. In this case the compatibility 
conditions \eq{compat} read
\[\big[\del_m \big(Z^j_{jk}\big)-\del_k\big(Z^j_{jm}\big)\big](h^\alpha-\kappa^j)\equiv0\,.\]
We claim that the first factor vanishes, i.e., \eq{R0Zj=i1} is  valid in this case as well. 
Indeed, the derivation of \eq{R0Zj=i} is valid in this case as well. By assumption
$\exists \alpha$ such that $\{m,k\}\subset A_\alpha$, such that 
$Z^m_{mj}=Z^k_{kj}=Z^\alpha_j$ due to assumption \eq{Zassum1}. Thus the compatibility 
conditions for \eq{ki} are met.

\bigskip

\emph{Compatibility conditions for \eq{hi}:}
For $\alpha=1,\dots,s_0$ we need to verify that
\beq\label{mixed-h}
	\del_k\del_m\ha^\alpha=\del_m\del_k \ha^\alpha\qquad \mbox{for $m\neq k$.}
\eeq 
For a fixed $\alpha$ there are several cases depending on whether both, exactly one, or neither 
of $m$ and $k$ belong to $A_\alpha$.

{\em $\bullet$ Case 1.} If both $m$ and $k$ belong to $A_\alpha$ then \eq{mixed-h} amounts to 
the trivial condition $0=0$.

{\em $\bullet$ Case 2.} Without loss of generality we assume that $m\in A_\alpha\not \ni k$.
The compatibility condition \eq{mixed-h} then reduces to 
\beq\label{case2}
	0 =\del_m\del_k \ha^\alpha=\del_m\big(W^\alpha_k(\ka^k-\ha^\alpha)\big)
	=\big[\del_m(W^\alpha_k)-W^\alpha_kZ^k_{km}\big](\ka^k-\ha^\alpha)\,.
\eeq
Since $|A_\alpha|>1$ there is a $j\neq m\in A_\alpha$, such that $W^\alpha_k=Z^j_{jk}$. 
We will verify \eq{case2} by showing that the first factor on the right-hand side vanishes:
\beq\label{van_rhs}
	\del_m(W^\alpha_k)-W^\alpha_kZ^k_{km}=\del_m(Z^j_{jk})-Z^j_{jk}Z^k_{km} = 0\,.
\eeq
Observing that \eq{Zassum2} again yields \eq{Zs}, we see that \eq{flat} gives
\bea
	\del_m\big(Z^j_{jk}\big) &=& Z^j_{jk}Z^k_{km}+Z^j_{jm}Z^m_{mk} - Z^j_{jm}Z^j_{jk}\nn\\
	&=& Z^j_{jk}Z^k_{km} + Z^j_{jm} (Z^m_{mk} - Z^j_{jk}) = Z^j_{jk}Z^k_{km}\,,\nn
\eea
where we have used that $Z^m_{mk} = Z^j_{jk}=W^\alpha_k$, since $j,\, m\in A_\alpha$.
This verifies \eq{van_rhs}.

{\em $\bullet$ Case 3.} If neither $m$ nor $k$ belong to $A_\alpha$ then the compatibility 
condition \eq{mixed-h} reduces to
\beq\label{last}
	\del_m\big[W^\alpha_k(\kappa^k-h^\alpha) \big]
	=\del_k\big[W^\alpha_m(\kappa^m-h^\alpha) \big]\,,
\eeq
where $W^\alpha_k=Z^j_{jk}$ and $W^\alpha_m=Z^j_{jm}$, for all $j\in A_\alpha$, and
with the provision that 
\[\kappa^k=h^\beta\qquad\mbox{if $k\in A_\beta$ ($\beta\neq \alpha$)} \qquad\mbox{and}\qquad 
\kappa^m=h^\gamma\qquad\mbox{if $m\in A_\gamma$ ($\gamma\neq \alpha$).}\]
There are three sub-cases depending on whether both, exactly one, or neither of  $m$ nor 
$k$ belong to some $A_\beta$ with $\beta\neq\alpha$. To verify the identity \eq{last} we
use the same techniques as above: flatness \eq{R0Z}, and \eq{Zassum1}-\eq{Zassum3}.
The details are similar to the previous cases and are left out.
\foorp

\noindent
Combining Lemmas \ref{rich_cnstrnd1} - \ref{int-hk-system} with Theorem \ref{dar3}
completes the proof of Proposition \ref{rich_cnstrnd}.

%
%
\section{Examples}\label{examples}

\subsection{The Euler system for 1-dimensional compressible 
flow}
\begin{example}\label{gas_dyn}	
The Euler system for 1-dimensional compressible 
flow in Lagrangian variables is:
\bea
	v_t - u_x &=& 0  \label{mass}\\
	u_t + p_x&=& 0  \label{mom}\\
	E_t + (up)_x &=& 0\,,\label{energy}
\eea
where $v,\, u,\, p$ are the specific volume, velocity, and pressure, respectively, 
and $E=\epsilon+\frac{u^2}{2}$ is the total specific energy, where $\epsilon$ denotes 
the specific internal energy.
The thermodynamic variables are related through Gibbs' relation $d\epsilon=TdS-pdv$,
where $T$ is temperature and $S$ is specific entropy. Using this, the Euler system
may be rewritten as
\bea
	v_t - u_x &=& 0 \label{Cont}\\
	u_t + p_x&=& 0  \label{Mom}\\
	S_t &=& 0\,. \label{Ener}
\eea
The system is closed by prescribing a pressure function $p=p(v,S)>0$, and we 
make the standard assumption that $p_v(v,S)<0$.
For a given a pressure function the eigenvalues of \eq{Cont}-\eq{Ener} are:
\[\lambda^1=-\sqrt{-p_v}\,,\qquad 
\lambda^2\equiv 0\,,\qquad 
\lambda^3=\sqrt{-p_v}\,,\]
with corresponding right and left eigenvectors (in $(v,u,S)$-space and normalized according to $R_i\cdot L^j\equiv \delta_i^j$):
\beq\label{revs}
	R_1=\left[\,1,\, \sqrt{-p_v},\, 0\,\right]^T\,,\qquad
	R_2=\left[\,-p_S,\, 0,\, p_v\,\right]^T\,,\qquad
	R_3=\left[\,1,\, -\sqrt{-p_v},\, 0\,\right]^T\,,
\eeq
and 
\beq\label{levs}
	L^1=\textstyle\frac{1}{2}\Big[\, 1\,,\, \frac{1}{\sqrt{-p_v}}\,,\, \frac{p_S}{p_v}\, \Big]\,,\qquad
	L^2=\Big[\, 0\,,\, 0\,,\, \frac{1}{p_v}\, \Big]\,,\qquad
	L^3=\textstyle\frac{1}{2}\Big[\, 1\,,\, -\frac{1}{\sqrt{-p_v}}\,,\, \frac{p_S}{p_v}\, \Big]\,.
\eeq
The fact that the eigenfields depend only on the pressure allows us to pose the ``inverse" problem:
\begin{problem}\label{prob3} (Problem \ref{prob1} for the Euler system)  For a given pressure 
function $p=p(v,S)>0$, with $p_v<0$, consider the frame given by \eq{revs}.
Then: determine the class of conservative systems with these as eigenfields
by solving the associated $\lambda$-system \eq{sev1}-\eq{sev2}
for the eigenvalues $\lambda^1,\, \lambda^2,\, \lambda^3$.
\end{problem}
\begin{remark}
	As explained in Remark \ref{main_rmk} the two forms \eq{mass}-\eq{energy} 
	and \eq{Cont}-\eq{Ener} of the Euler system are equivalent as far as Problem \ref{prob3} 
	is concerned. However, they are not equivalent at the level of computing weak solutions to 
	Cauchy problems \cite{sm}.
\end{remark}
We shall see that there are two distinct cases depending on whether the frame is rich or not.
From \eq{revs}, \eq{levs}, and \eq{christoffel} we obtain
\[\Gamma_{21}^2 = \Gamma_{23}^2 = 0\,,\qquad 
\Gamma_{31}^3 = \Gamma_{13}^1 = \textstyle-\frac{p_{vv}}{4p_{v}}\,,\qquad
\Gamma_{12}^1 = \Gamma_{32}^3 = -\frac{p_v}{2} \big(\frac{p_S}{p_v}\big)_v\,.\]
Thus the PDEs \eq{sev1} in the $\lambda$-system are given by
\bea
	r_1(\lambda^2) &=& 0 \label{a}\\
	r_1(\lambda^3) &=& \textstyle\frac{p_{vv}}{4p_{v}}(\lambda^3-\lambda^1) \label{b} \\
	r_2(\lambda^1) &=& \textstyle\frac{p_v}{2} \big(\frac{p_S}{p_v}\big)_v(\lambda^1-\lambda^2)
	\label{c}\\
	r_2(\lambda^3) &=& \textstyle\frac{p_v}{2} \big(\frac{p_S}{p_v}\big)_v(\lambda^3-\lambda^2)
	\label{d}\\
	r_3(\lambda^1) &=& \textstyle\frac{p_{vv}}{4p_{v}}(\lambda^1-\lambda^3) \label{e} \\
	r_3(\lambda^2) &=& 0 \,. \label{f}
\eea
where $r_1,\, r_2,\, r_3$ are the derivations corresponding to $R_1,\, R_2,\, R_3$.
The coefficients appearing in the algebraic equations \eq{sev2} of the $\lambda$-system are:
\[\Gamma_{32}^1 = \Gamma_{12}^3 = \textstyle -\frac{p_v}{2} \big(\frac{p_S}{p_v}\big)_v\,,\qquad 
\Gamma_{23}^1 = \Gamma_{21}^3 = -\frac{p_v}{4} \big(\frac{p_S}{p_v}\big)_v\,,\qquad 
\Gamma_{31}^2 = \Gamma_{13}^2 = 0\,. \]
Using these we see that the matrix $N$ in \eq{3x3_alg_cond} is
given by
\beq\label{N}
	N=\textstyle \frac{p_v}{4} \big(\frac{p_S}{p_v}\big)_v
	\left[\begin{array}{rr}
	-2 & 1\\
	0 & 0 \\
	-2 & 1\end{array}\right]\,. 
\eeq
By assumption $p_v < 0$, and we thus consider the two cases: 
\begin{itemize}
\item [$(a)$] $\big(\frac{p_S}{p_v}\big)_v\equiv 0$: $\rank(N)=0$ and there are no algebraic
constraints among the eigenvalues. Proposition \ref{gen_case} and Theorem \ref{rich_unconstr}  
show that the Euler system is rich in this case, and that the general solution 
to Problem \ref{prob3} depends on three functions of one variable.
\item [$(b)$] $\big(\frac{p_S}{p_v}\big)_v\neq 0$: $\rank(N)=1$ and there 
is a single algebraic relation among the eigenvalues: 
\beq\label{alg_rel}
	\lambda^1+\lambda^3=2\lambda^2\,.
\eeq
\end{itemize}
For case (a) the pressure must be of the form $p(v,S)=\Pi(\xi)$, where $\xi=v+F(S)$ and $F$ is 
a function of $S$ alone. It follows from \eq{a} and \eq{f} that $\lam^2$ is any function 
of $S$ alone, and from \eq{c} and \eq{d} it follows that $\lam^1=A(\xi,u)$ and 
$\lam^3=B(\xi,u)$. The functions $A$ and $B$ should then be determined from
the remaining equations  \eq{b} and \eq{e}: 
\[A_\xi-\sqrt{-\Pi'(\xi)}A_u=\alpha(B-A)\,,\qquad B_\xi+\sqrt{-\Pi'(\xi)}B_u=\alpha(A-B)\,,\]
where $\alpha=-\textstyle\frac{p_{vv}}{4p_{v}}=-\textstyle\frac{\Pi''(\xi)}{4\Pi'(\xi)}$. These 
equations must be solved case-by-case. However, by deriving a second order equation for 
$A$ alone, which in turn will determine $B$, we see that solving for $A$ and $B$ requires
two functions of one variable. Together with $\lam^2(S)$ these are the three functions that 
determine a general solution in case (a).

We proceed to explicitly solve the PDE system in the more interesting second case (b).
Adding \eq{c} and \eq{d}, and using \eq{alg_rel}, gives $r_2(\lambda^1+\lambda^3) = 0$.
Applying \eq{alg_rel} again yields $r_2(\lambda^2) = 0$.
Together with \eq{a} and \eq{f} this shows that
\beq\label{2=const}
	\lambda^2\equiv \bar \lambda \quad \mbox{(constant)}.
\eeq
To continue it is convenient to express the coordinate frame $\{\del_v,\, \del_u,\, \del_S\}$ 
in terms of the $r$-frame:
\beq\label{coord_frame}
	\partial_v = \textstyle\frac 1 2(r_1+r_3) \,,\qquad
	\partial_u = \frac 1 {2\sqrt{-p_v}}(r_1-r_3) \,,\qquad
	\partial_S = \frac 1 {p_v}r_2+\frac {p_S} {2 p_v}(r_1+r_3) \,.
\eeq
A direct calculation, using \eq{coord_frame} and \eq{alg_rel},
now shows that $\Delta:=(\lambda^3-\lambda^1)$ satisfies $\del_u\Delta = 0$ and
\[\del_v\Delta = \textstyle\big(\frac{p_{vv}}{2p_v}\big)\Delta \,,\qquad
\del_S\Delta = \big(\frac{p_{vS}}{2p_v}\big) \Delta\,.\]
Integration with respect to $v$ and $S$, respectively, shows that 
\[\Delta(v,S)=\lambda^3-\lambda^1=D\sqrt{-p_v}\]
for a constant $D$. Invoking \eq{alg_rel} a last time (with $\lambda^2\equiv \bar\lambda$), 
we conclude that the general solution of the $\lambda$-system for non-rich 
gas dynamics (case (b)) depends on two constants $\bar\lam$, $C$ according to:
\[\lambda^1 = \bar \lambda-C\sqrt{-p_v}  \,,\qquad
\lambda^2 \equiv \bar \lambda \,,\qquad
\lambda^3 = \bar \lambda+C\sqrt{-p_v}\,. \]
\end{example}

\subsection{Additional examples with non-rich frames }\label{ex_non-rich}
The case of non-rich gas dynamics in Example \ref{gas_dyn} provides an example of 
Subcase IIa in Section \ref{rank1_n=3} in which the compatibility conditions 
\eq{frob-comp1}-\eq{frob-comp2} for the Frobenius system \eq{IIa-solved} hold identically.
The following Example~\ref{nonrich_IIa_trivial}  of a system in the same category, but with 
only trivial solutions, shows that these compatibility conditions do not necessarily hold as 
identities and must be checked in each case. This verifies the last claim in Proposition 
\ref{case_IIa}.

\begin{example}\label{nonrich_IIa_trivial} 
{\sc$n= 3$, non-rich system with a single algebraic relation of type} $\mathrm{IIa}$ {\sc and only  trivial solutions.}
We prescribe the frame $\{R_1,\, R_2,\, R_3\}$ on $\RR^3$ by 
\[R_1=[0,0,1]^T\,,\qquad R_2=[0,1,u^1]^T\,,\qquad R_3=[u^3,0,1]^T\,.\] 
The system \eq{pde1}-\eq{pde6} takes form
\begin{eqnarray*}
\del_3\lambda^2&=&0,\\
\del_3\lambda^3&=&0,\\
\del_2\lambda^1+u^1\del_3\lambda^1&=&0,\\
\del_2\lambda^3+u^1\del_3\lambda^3&=&0,\\
u^3\del_1\lambda^1+\del_3\lam^1+\textstyle\frac 1 u^3(\lam^3-\lam^1)&=&0,\\
u^3\del_1\lambda^2+\del_3\lam^2&=&0.
\end{eqnarray*}
In this case the matrix  $N$ is given by 
\[N=
	\left[\begin{array}{cc}
	u^3&\frac {u^1 }{u^3}  \\
	0  &0\\
	0&0
	\end{array}\right]\]
and has rank 1.
It corresponds to the unique algebraic relation that involves all three eigenvalues:
\[(u^3)^2(\lam^2-\lam^1)+u^1(\lam^3-\lam^1)=0.\]
The only solution of the above  system of differential and algebraic conditions is trivial 
$\lam^1=\lam^2=\lam^3=const.$  
\end{example} 

Next we provide examples of a non-rich 3-frame with a single algebraic relation of type IIb.
The examples show that in accordance with Proposition~\ref{case_IIb} the solution of the 
$\lam$-system in this case either depends on one arbitrary function  of one variable and one arbitrary constant, or is trivial. 
\begin{example}\label{nonrich_IIb_non-trivial}{\sc$n= 3$, non-rich system with a single algebraic 
relation of type} $\mathrm{IIb}$ {\sc and non trivial solutions}.
Consider a system with eigenvectors:
$$R_1=(0,1,0)^T,\quad  R_2=(1,0,0)^T, \quad R_3=(u^2,u^3,1)^T.$$ 
Since $[r_1,r_3]=r_2$, 
these vector-fields are not pair-wise 
in involution, and therefore the system is not rich. 
The $\lam$-system implies a unique algebraic relations $\lam^2=\lam^3$, which involves only
two of the $\lam$'s. We are thus in Subcase IIb in Section \ref{subcaseIIb}. A computation show 
that  $\Gamma^3_{31}=0$ and $\Gamma^2_{21}=0$ are equal and thus according to 
Proposition~\ref{case_IIb} we can expect a non-trivial general solution depending on one arbitrary 
function  of one variable and one arbitrary constant. 
Indeed, taking into account the algebraic relation the differential equations become:  
\begin{eqnarray*}
r_2(\lambda^1)&=&\del_1\lambda^1=0,\\
r_3(\lambda^1)&=&u^2\,\del_1\lambda^1+u^3\,\del_2\lambda^1+\del_3\lam^1=0,\\
r_1(\lambda^2)&=&\del_2\lambda^2=0,\\
r_2(\lambda^2)&=&\del_1\lambda^2=0,\\
r_3(\lambda^2)&=&u^2\,\del_1\lambda^2+u^3\,\del_2\lambda^2+\del_3\lam^2=0.
\end{eqnarray*}
The general solution of this system is $\lam^1=\varphi(u_3^2-2\,u_2)\,$,  $\lambda^2=\lambda^3\equiv \bar\lambda$, 
where   $\bar\lambda$ is an arbitrary constant and 
$\varphi$ is an arbitrary function of one variable.
\end{example}
\begin{example}\label{nonrich_IIb_trivial}{\sc$n= 3$, non-rich system with a single algebraic 
relation of type} $\mathrm{IIb}$ {\sc and only trivial solutions}.
Consider a system with eigenvectors:
$$R_1=(u^1+u^3,1,0)^T,\quad  R_2=(1,0,0)^T, \quad R_3=(u^2,u^3,-u^2)^T.$$ 
Since $[r_1,r_3]=-\frac{u^3}{u^2}\,r^1+\frac{(u^1+u^3)\,u^3}{u^2}\,r^2+\frac{1}{u^2}\,r^3$ these vector-fields are not pair-wise 
in involution, and therefore the system is not rich. 
The $\lam$-system implies a unique algebraic relations $\lam^2=\lam^3$, which involves only
two of the $\lam$'s. We are thus in Subcase IIb in Section \ref{subcaseIIb}. A computation shows that  
$\Gamma^3_{31}=0$ and $\Gamma^2_{21}=1$ are not equal and thus according to Proposition~\ref{case_IIb} 
we can expect only trivial general solution.
Indeed, taking into account the algebraic relation $\lam^2=\lam^3$ the differential equations
$r_1(\lam^2)= \Gamma^2_{21}(\lam^1-\lam^2)$ and  $r_1(\lam^3)= \Gamma^3_{31}(\lam^1-\lam^3)$ imply that   
\beq
r_1(\lambda^2)=(\lam^2-\lam^1)\quad\mbox{and}\quad r_1(\lambda^2)=0,\eeq
which implies $\lam^1=\lam^2$. The only non-trivial solution of the $\lam$ system is therefore  $\lam^1=\lam^2=\lam^3=\bar \lam$,
where   $\bar\lambda$ is an arbitrary constant.

\end{example}
We include  an example which illustrates the fact that maximal rank of $N$ ($\rank(N)=n-1$)
implies that all solutions of the $\lambda$-system are trivial (see section \ref{rank_discuss}).
\begin{example}\label{nonrich_n4_maximalrank}{\sc $n=4$, system with maximal rank of algebraic constraints ($\rank N=3)$.}
The $\lambda$-system for the following non-rich frame in $\RR^4$  
$$R _1=(1, 0, u^2, u^4), \quad R_2=(0, 1, u^1, 0), \quad R_3=(u^3, 0, 1, 0),  \quad R_4=(1, 0, 0, 0)$$
includes three independent algebraic relations:
$$\lam^2=\lam^1,\quad  \lam^3=\lam^4,\quad \lam^3=\lam^2\,,$$
which immediately shows that all solutions must be trivial:
$\lam^1=\lam^2=\lam^3=\lam^4\equiv const.$
\end{example}

\subsection{Additional examples with rich frames}\label{ex_rich}
\begin{example}\label{n2} {\sc $n=2$ system.}
Suppose we want to find the most general system \eq{claw} of two conservation laws whose
eigencurves are hyperbolas  and radial straight lines. As a $2\times 2$-system 
it is necessarily rich and as Riemann invariants we use 
$w^1:=u^1u^2$ and $w^2:=\frac {u_2} {u_1}$, where $(u^1,u^2)$ are standard   coordinates on $\RR^2$. 
Then the matrix of left eigenvectors is 
\[ L(u^1,u^2)=\left[\begin{array}{c}
    \nabla w^1\\
    \nabla w^2
    \end{array}\right]\, =  \left[\begin{array}{cc}
   u^2 & u^1\\
    -\frac {u^2} {(u^1)^2} & \frac 1 {u^1}
    \end{array}\right]\,\]
and the matrix of right  eigenvectors is $R=L^{-1}$.
A direct computation shows that  $Z^2_{21}=\frac 1 {2w^1}$ and $Z^1_{12}=0$, 
such that the $\lambda$-system \eq{sev1rich} becomes
\begin{eqnarray}
	\label{k1} \del_2 \ka^1(w^1,w^2)  & = &  0\\
	\label{k2} \del_1\ka^2(w^1,w^2)  & = &  \frac 1 {2w^1}(\ka^1-\ka^2).
\end{eqnarray}
By \eq{k1} the first eigenvalue is an arbitrary function of $w^1$ alone, $\ka^1=\ka^1(w^1)$, 
and \eq{k2} then yields
\[\ka^2(w^1,w^2)=\frac{1}{\sqrt {w^1}}\big[h(w^2)+g(w^1)\big]\,,\qquad
\mbox{where}\quad g(w^1)=\int_1^{w^1}\frac{\ka^1(\xi)}{\sqrt{\xi}}\,d\xi,\]
and $h$ is an arbitrary function of $w^2$ alone. In accordance with Theorem 
\ref{rich_unconstr} we see that the general solution depends on two functions of one 
variable.
\end{example}

As we proved in Section~\ref{rich-orthogonal}, for a rich and orthogonal frame the algebraic 
part of the $\lam$-system becomes trivial.  We consider a concrete case in this class.

\begin{example}\label{rich_orth} {\sc rich orthogonal frame.} Suppose we want to find the most 
general system \eq{claw} of two conservation laws whose
eigencurves are the coordinates curves of spherical coordinates (radial lines, latitudes, and longitudes). 
The corresponding Riemann invariants are polar coordinates $(r, \theta,\phi)$. More precisely  we let
$\rho: \mathbb R^3 \setminus \{(u^1,0,u^3)\, |\, u^1\leq 0,\, z\leq 0\,\} \to
\mathbb R_+\times (-\pi,\pi)\times \mathbb (0,\pi)$ be the change
from rectangular Cartesian coordinates $(u^1,u^2,u^2)$ to spherical polar
coordinates $(r,\theta,\phi)$,
\beq\label{chng_coords_2}
    \rho(u,v,w)=(r,\theta,\phi),
\eeq
where
\[  r =\sqrt{(u^1)^2+(u^2)^2+(u^3)^2},\quad
    \theta = \arctan\left(\frac{\sqrt{(u^1)^2+(u^2)^2}}{u^3}\right),\quad\mbox{and}\quad
    \phi = \arctan\left(\frac {u^2}{ u^1}\right).\]
Defining $R$ by $R^{-1}=D\rho$ we have
\[R=\left[\begin{array}{ccc}
\sin\theta\cos\phi & r\cos\theta \cos\phi & -r\sin\theta \sin\phi \\
\sin\theta\sin\phi & r\cos\theta \sin\phi  & r\sin\theta\cos\phi\\
\cos\theta & -r\sin\theta & 0
\end{array}\right]\,.\]

A calculation now yields that $Z^1_{12}=Z^1_{13}=Z^2_{23}=0$, $ Z^2_{21}=Z^3_{31}=\frac 1 r$, 
$Z^{3}_{32}=\arctan(\theta)$  and, as expected,  all $Z$'s with three distinct indices are zero. Thus
according to \eq{sev1rich}-\eq{sev2rich}, the $\lambda$-system contains only the following set of differential equations:
\begin{eqnarray*}
    \del_\theta \ka^1=0\hskip16mm& \quad &\del_\phi \ka^1=0\, ,\\
\del_r  \ka^2=\frac 1 r (\ka^1-\ka^2)\, ,& \quad &\del_\phi \ka^2=0\, ,\\
\del_r \ka^3=\frac 1 r (\ka^1-\ka^3)\, ,& \quad &\del_\theta \ka^3=\arctan\theta(\ka^2-\ka^3).
\end{eqnarray*}
The general solution of this system depends on three arbitrary functions $F, G, H$ each depending on a single variable: 
\begin{eqnarray*}
\ka^1&=&F(r)\, ,\\
\ka^2&=&\frac{H(\theta)+g(r)}{r}, \text{ where } g(r):=\int_1^r F(\xi)\, d\xi\, ,\\
\ka^3& =& \frac{K(\phi) + \int_{\theta_0}^\theta
    H(\xi)\cos \xi \, d\xi}{r\sin\theta}+\frac{g(r)}{r}.
\end{eqnarray*}
\end{example}

\begin{example}\label{constant_eigenfields} 
{\sc Constant eigenfields.} If the frame $\{R_1(u),\dots,R_n(u)\}$ consists of constant 
vector fields $R_1(u)\equiv R_i$ then it is rich according to Definition \ref{rich_sys}. 
Furthermore, all connection coefficients $Z_{ij}^k$ vanish and Theorem \ref{rich_unconstr} 
applies. As associated Riemann invariants we may choose $w^i(u)=L^i\cdot u$, where each 
$L_i$ is constant. 
A calculation shows that the general solution of $\lambda$-system \eq{sev1rich} 
is given by $\ka^1=\phi^1(w^1),\dots,\phi^n(w^n)$, where $\phi^1,\dots,\phi^n$ 
are arbitrary functions of one variable, in accordance with Theorem \ref{rich_unconstr}. 
That is, the eigenvalues solving Problem 
\ref{prob1} in this case are given by $\lam^i(u)=\phi^i(L^i\cdot u)$, $i=1,\dots,n$.
\end{example}

We have seen an example of  rich  $3\times 3$-system with no algebraic constraints ($\rank N=0$)  
when we analyzed eigenvector-fields of the  Euler system in  Example~\ref{gas_dyn}. The following two examples
illustrate our observation  that Richness $\not\Rightarrow \rank(N)=0$. In Example~\ref{rich_rank1} we have a rich system
with $\rank(N)=1$ with non trivial general solution depending on 1 arbitrary constant and 1 arbitrary  function of 1 variable.
In Example~\ref{rich_rank2} we have a rich system with $\rank(N)=2$ and therefore only  trivial solutions.

\begin{example}\label{rich_rank1} {\sc $n=3$, rich, $\rank (N)=1$, non-trivial.}
Consider a system with eigenvectors:
$$R_1=(1,0,u^2)^T, \quad R_2=(0,1,u^1)^T, \quad R_3=(0,0,-1)^T.$$ 
These vectors commute and therefore any system \eq{claw} with these as eigenvectors is rich. 
The only non-zero connection components are $\Gamma_{12}^3=\Gamma_{21}^3=-1$,
such that the $\lam$-system implies a single algebraic relation: $\lam^1=\lam^2$. 
Thus this example illustrates the results in Section \ref{rich-rank-non-zero}. 
Taking into account this relation the differential equations become:  
\begin{eqnarray*}
r_1(\lambda^1)&=&\del_1\lambda^1+u^2\del_3\lambda^1=0,\\
r_1(\lambda^3)&=&\del_1\lambda^3+u^2\del_3\lambda^3=0,\\
r_2(\lambda^1)&=&\del_2\lambda^1+u^1\del_3\lambda^1=0,\\
r_2(\lambda^3)&=&\del_2\lambda^3+u^1\del_3\lambda^3=0,\\
r_3(\lambda^1)&=&-\del_3\lambda^1=0.
\end{eqnarray*}
The general solution of this system is  $\lambda^1=\lambda^2\equiv \bar\lambda$, and 
$\lambda^3=\varphi(u^3-u^1u^2)$, where $\bar\lambda$ is an arbitrary constant and 
$\varphi$ is an arbitrary function of one variable. 
\end{example}

\begin{example}\label{rich_rank2} {\sc $n=3$, rich, $\rank N=2$, trivial.}
Consider the vector fields
\[R_1=[u^1,\, u^2,\, 0]^T\,,\qquad R_2= [-u^2,\, u^1,\, 0]^T\,,\qquad R_3=[-u^2,\, u^1,\, 1]^T\,.\]
We are then searching for systems \eq{claw} whose eigencurves are 
horizontal lines through the $u^3$-axis, horizontal circles, and vertical helices.
A calculation shows that these determine a rich frame. 
We can rewrite this problem in terms of Riemann invariants using the change of coordinates: 
$\rho:\RR^3\to\RR^3$ given by
\beq\label{chng_coords_6}
    \rho^{-1}(r,\theta,h)=(u^1,u^2,u^3)
    =\big(r\cos(\theta+h),r\sin(\theta+h),h\big)\,.
\eeq
A calculation show  that $Z^1_{23}=-r$ and 
$Z^2_{13}=\frac 1 r$. It follows from \eq{sev2rich} 
and Proposition \ref{equal_lambdas} that the only solutions are the trivial 
solutions $\lam^1=\lam^2=\lam^3\equiv const.$
\end{example}



\begin{bibdiv}
\begin{biblist}
\bib{bj}{article}{
   author={Baiti, Paolo},
   author={Jenssen, Helge Kristian},
   title={Blowup in $L\sp \infty$ for a class of genuinely nonlinear
   hyperbolic systems of conservation laws},
   journal={Discrete Contin. Dynam. Systems},
   volume={7},
   date={2001},
   number={4},
   pages={837--853},
   issn={1078-0947},
   review={\MR{1849664 (2003m:35155)}},
}
\bib{bi1}{article}{
   author={Bianchini, Stefano},
   title={BV solutions of the semidiscrete upwind scheme},
   journal={Arch. Ration. Mech. Anal.},
   volume={167},
   date={2003},
   number={1},
   pages={1--81},
   issn={0003-9527},
   review={\MR{1967667 (2004k:35249)}},
}\bib{bi2}{article}{
   author={Bianchini, Stefano},
   title={Hyperbolic limit of the Jin-Xin relaxation model},
   journal={Comm. Pure Appl. Math.},
   volume={59},
   date={2006},
   number={5},
   pages={688--753},
   issn={0010-3640},
   review={\MR{2172805}},
}
\bib{bb}{article}{
   author={Bianchini, Stefano},
   author={Bressan, Alberto},
   title={Vanishing viscosity solutions of nonlinear hyperbolic systems},
   journal={Ann. of Math. (2)},
   volume={161},
   date={2005},
   number={1},
   pages={223--342},
   issn={0003-486X},
   review={\MR{2150387 (2007i:35160)}},
}
\bib{br1}{article}{
   author={Bressan, Alberto},
   title={Global solutions of systems of conservation laws by wave-front
   tracking},
   journal={J. Math. Anal. Appl.},
   volume={170},
   date={1992},
   number={2},
   pages={414--432},
   issn={0022-247X},
   review={\MR{1188562 (93k:35166)}},
}
\bib{br}{book}{
   author={Bressan, Alberto},
   title={Hyperbolic systems of conservation laws},
   series={Oxford Lecture Series in Mathematics and its Applications},
   volume={20},
   note={The one-dimensional Cauchy problem},
   publisher={Oxford University Press},
   place={Oxford},
   date={2000},
   pages={xii+250},
   isbn={0-19-850700-3},
   review={\MR{1816648 (2002d:35002)}},
}
\bib{bcggg}{book}{
   author={Bryant, R. L.},
   author={Chern, S. S.},
   author={Gardner, R. B.},
   author={Goldschmidt, H. L.},
   author={Griffiths, P. A.},
   title={Exterior differential systems},
   series={Mathematical Sciences Research Institute Publications},
   volume={18},
   publisher={Springer-Verlag},
   place={New York},
   date={1991},
   pages={viii+475},
   isbn={0-387-97411-3},
   review={\MR{1083148 (92h:58007)}},
}
\bib{chen}{article}{
   author={Chen, Gui-Qiang},
   title={Compactness methods and nonlinear hyperbolic conservation laws},
   conference={
      title={Some current topics on nonlinear conservation laws},
   },
   book={
      series={AMS/IP Stud. Adv. Math.},
      volume={15},
      publisher={Amer. Math. Soc.},
      place={Providence, RI},
   },
   date={2000},
   pages={33--75},
   review={\MR{1767623 (2001f:35249)}},
}
\bib{amf}{book}{
   author={Choquet-Bruhat, Yvonne},
   author={DeWitt-Morette, C{\'e}cile},
   author={Dillard-Bleick, Margaret},
   title={Analysis, manifolds and physics},
   edition={2},
   publisher={North-Holland Publishing Co.},
   place={Amsterdam},
   date={1982},
   pages={xx+630},
   isbn={0-444-86017-7},
   review={\MR{685274 (84a:58002)}},
}
\bib{cl}{article}{
   author={Conlon, Joseph G.},
   author={Liu, Tai Ping},
   title={Admissibility criteria for hyperbolic conservation laws},
   journal={Indiana Univ. Math. J.},
   volume={30},
   date={1981},
   number={5},
   pages={641--652},
   issn={0022-2518},
   review={\MR{625595 (82j:35095)}},
}
\bib{daf1}{article}{
   author={Dafermos, C. M.},
   title={Stability for systems of conservation laws in several space
   dimensions},
   journal={SIAM J. Math. Anal.},
   volume={26},
   date={1995},
   number={6},
   pages={1403--1414},
   issn={0036-1410},
   review={\MR{1356450 (96i:35079)}},
}
\bib{daf}{book}{
   author={Dafermos, Constantine M.},
   title={Hyperbolic conservation laws in continuum physics},
   series={Grundlehren der Mathematischen Wissenschaften [Fundamental
   Principles of Mathematical Sciences]},
   volume={325},
   edition={2},
   publisher={Springer-Verlag},
   place={Berlin},
   date={2005},
   pages={xx+626},
   isbn={978-3-540-25452-2},
   isbn={3-540-25452-8},
   review={\MR{2169977 (2006d:35159)}},
}
\bib{dar}{book}{
   author={Darboux, Gaston},
   title={Le\c cons sur les syst\`emes orthogonaux et les coordonn\'ees
   curvilignes. Principes de g\'eom\'etrie analytique},
   language={French},
   series={Les Grands Classiques Gauthier-Villars. [Gauthier-Villars Great
   Classics]},
   note={The first title is a reprint of the second (1910) edition; the
   second title is a reprint of the 1917 original;
   Cours de G\'eom\'etrie de la Facult\'e des Sciences. [Course on Geometry
   of the Faculty of Science]},
   publisher={\'Editions Jacques Gabay},
   place={Sceaux},
   date={1993},
   pages={600},
   isbn={2-87647-016-0},
   review={\MR{1365963 (97b:01026)}},
}
\bib{dip1}{article}{
   author={DiPerna, Ronald J.},
   title={Global existence of solutions to nonlinear hyperbolic systems of
   conservation laws},
   journal={J. Differential Equations},
   volume={20},
   date={1976},
   number={1},
   pages={187--212},
   issn={0022-0396},
   review={\MR{0404872 (53 \#8672)}},
}
\bib{dip2}{article}{
   author={DiPerna, Ronald J.},
   title={Compensated compactness and general systems of conservation laws},
   journal={Trans. Amer. Math. Soc.},
   volume={292},
   date={1985},
   number={2},
   pages={383--420},
   issn={0002-9947},
   review={\MR{808729 (87g:35148)}},
}
\bib{fl}{article}{
   author={Friedrichs, K. O.},
   author={Lax, P. D.},
   title={Systems of conservation equations with a convex extension},
   journal={Proc. Nat. Acad. Sci. U.S.A.},
   volume={68},
   date={1971},
   pages={1686--1688},
   review={\MR{0285799 (44 \#3016)}},
}
\bib{gl}{article}{
   author={Glimm, James},
   title={Solutions in the large for nonlinear hyperbolic systems of
   equations},
   journal={Comm. Pure Appl. Math.},
   volume={18},
   date={1965},
   pages={697--715},
   issn={0010-3640},
   review={\MR{0194770 (33 \#2976)}},
}
\bib{gla}{book}{
   author={Glimm, James},
   author={Lax, Peter D.},
   title={Decay of solutions of systems of nonlinear hyperbolic conservation
   laws},
   series={Memoirs of the American Mathematical Society, No. 101},
   publisher={American Mathematical Society},
   place={Providence, R.I.},
   date={1970},
   pages={xvii+112},
   review={\MR{0265767 (42 \#676)}},
}
\bib{go}{article}{
   author={Godunov, S. K.},
   title={An interesting class of quasi-linear systems},
   language={Russian},
   journal={Dokl. Akad. Nauk SSSR},
   volume={139},
   date={1961},
   pages={521--523},
   issn={0002-3264},
   review={\MR{0131653 (24 \#A1501)}},
}
\bib{il}{book}{
   author={Ivey, Thomas A.},
   author={Landsberg, J. M.},
   title={Cartan for beginners: differential geometry via moving frames and
   exterior differential systems},
   series={Graduate Studies in Mathematics},
   volume={61},
   publisher={American Mathematical Society},
   place={Providence, RI},
   date={2003},
   pages={xiv+378},
   isbn={0-8218-3375-8},
   review={\MR{2003610 (2004g:53002)}},
}
\bib{je}{article}{
   author={Jenssen, Helge Kristian},
   title={Blowup for systems of conservation laws},
   journal={SIAM J. Math. Anal.},
   volume={31},
   date={2000},
   number={4},
   pages={894--908 (electronic)},
   issn={0036-1410},
   review={\MR{1752421 (2001a:35114)}},
}
\bib{jy}{article}{
   author={Jenssen, Helge Kristian},
   author={Young, Robin},
   title={Gradient driven and singular flux blowup of smooth solutions to
   hyperbolic systems of conservation laws},
   journal={J. Hyperbolic Differ. Equ.},
   volume={1},
   date={2004},
   number={4},
   pages={627--641},
   issn={0219-8916},
   review={\MR{2111577 (2006e:35225)}},
}
\bib{jo}{book}{
   author={John, Fritz},
   title={Partial differential equations},
   series={Applied Mathematical Sciences},
   volume={1},
   edition={4},
   publisher={Springer-Verlag},
   place={New York},
   date={1991},
   pages={x+249},
   isbn={0-387-90609-6},
   review={\MR{1185075 (93f:35001)}},
}
\bib{jmr}{article}{
   author={Joly, J.-L.},
   author={M{\'e}tivier, G.},
   author={Rauch, J.},
   title={A nonlinear instability for $3\times 3$ systems of conservation
   laws},
   journal={Comm. Math. Phys.},
   volume={162},
   date={1994},
   number={1},
   pages={47--59},
   issn={0010-3616},
   review={\MR{1272766 (95f:35145)}},
}
\bib{kk}{article}{
   author={Keyfitz, Barbara L.},
   author={Kranzer, Herbert C.},
   title={A system of nonstrictly hyperbolic conservation laws arising in
   elasticity theory},
   journal={Arch. Rational Mech. Anal.},
   volume={72},
   date={1979/80},
   number={3},
   pages={219--241},
   issn={0003-9527},
   review={\MR{549642 (80k:35050)}},
}
\bib{kru}{article}{
   author={Kru{\v{z}}kov, S. N.},
   title={First order quasilinear equations with several independent
   variables. },
   language={Russian},
   journal={Mat. Sb. (N.S.)},
   volume={81 (123)},
   date={1970},
   pages={228--255},
   review={\MR{0267257 (42 \#2159)}},
}
\bib{lax}{article}{
   author={Lax, P. D.},
   title={Hyperbolic systems of conservation laws. II},
   journal={Comm. Pure Appl. Math.},
   volume={10},
   date={1957},
   pages={537--566},
   issn={0010-3640},
   review={\MR{0093653 (20 \#176)}},
}
\bib{lee}{book}{
   author={Lee, John M.},
   title={Riemannian manifolds},
   series={Graduate Texts in Mathematics},
   volume={176},
   note={An introduction to curvature},
   publisher={Springer-Verlag},
   place={New York},
   date={1997},
   pages={xvi+224},
   isbn={0-387-98271-X},
   review={\MR{1468735 (98d:53001)}},
}
\bib{liu1}{article}{
   author={Liu, Tai Ping},
   title={The deterministic version of the Glimm scheme},
   journal={Comm. Math. Phys.},
   volume={57},
   date={1977},
   number={2},
   pages={135--148},
   issn={0010-3616},
   review={\MR{0470508 (57 \#10259)}},
}

\bib{ni}{article}{
   author={Nishida, Takaaki},
   title={Global solution for an initial boundary value problem of a
   quasilinear hyperbolic system},
   journal={Proc. Japan Acad.},
   volume={44},
   date={1968},
   pages={642--646},
   issn={0021-4280},
   review={\MR{0236526 (38 \#4821)}},
}

\bib{nism}{article}{
   author={Nishida, Takaaki},
   author={Smoller, Joel A.},
   title={Solutions in the large for some nonlinear hyperbolic conservation
   laws},
   journal={Comm. Pure Appl. Math.},
   volume={26},
   date={1973},
   pages={183--200},
   issn={0010-3640},
   review={\MR{0330789 (48 \#9126)}},
}

\bib{ol}{book}{
   author={Olver, Peter J.},
   title={Applications of Lie groups to differential equations},
   series={Graduate Texts in Mathematics},
   volume={107},
   edition={2},
   publisher={Springer-Verlag},
   place={New York},
   date={1993},
   pages={xxviii+513},
   isbn={0-387-94007-3},
   isbn={0-387-95000-1},
   review={\MR{1240056 (94g:58260)}},
}
\bib{ri}{article}{
   author={Risebro, Nils Henrik},
   title={A front-tracking alternative to the random choice method},
   journal={Proc. Amer. Math. Soc.},
   volume={117},
   date={1993},
   number={4},
   pages={1125--1139},
   issn={0002-9939},
   review={\MR{1120511 (93e:35071)}},
}
\bib{serre}{book}{
   author={Serre, Denis},
   title={Systems of conservation laws. 1},
   note={Hyperbolicity, entropies, shock waves;
   Translated from the 1996 French original by I. N. Sneddon},
   publisher={Cambridge University Press},
   place={Cambridge},
   date={1999},
   pages={xxii+263},
   isbn={0-521-58233-4},
   review={\MR{1707279 (2000g:35142)}},
}
\bib{serre2}{book}{
   author={Serre, Denis},
   title={Systems of conservation laws. 2},
   note={Geometric structures, oscillations, and initial-boundary value
   problems;
   Translated from the 1996 French original by I. N. Sneddon},
   publisher={Cambridge University Press},
   place={Cambridge},
   date={2000},
   pages={xii+269},
   isbn={0-521-63330-3},
   review={\MR{1775057 (2001c:35146)}},
}
\bib{sev}{article}{
   author={S{\'e}vennec, Bruno},
   title={G\'eom\'etrie des syst\`emes hyperboliques de lois de
   conservation},
   language={French, with English and French summaries},
   journal={M\'em. Soc. Math. France (N.S.)},
   number={56},
   date={1994},
   pages={125},
   issn={0037-9484},
   review={\MR{1259465 (95g:35123)}},
}

\bib{sever}{article}{
   author={Sever, Michael},
   title={Distribution solutions of nonlinear systems of conservation laws},
   journal={Mem. Amer. Math. Soc.},
   volume={190},
   date={2007},
   number={889},
   pages={viii+163},
   issn={0065-9266},
   review={\MR{2355635}},
}
\bib{sm}{book}{
   author={Smoller, Joel},
   title={Shock waves and reaction-diffusion equations},
   series={Grundlehren der Mathematischen Wissenschaften [Fundamental
   Principles of Mathematical Sciences]},
   volume={258},
   edition={2},
   publisher={Springer-Verlag},
   place={New York},
   date={1994},
   pages={xxiv+632},
   isbn={0-387-94259-9},
   review={\MR{1301779 (95g:35002)}},
}
\bib{sp}{book}{
   author={Spivak, Michael},
   title={A comprehensive introduction to differential geometry. Vol. I},
   edition={2},
   publisher={Publish or Perish Inc.},
   place={Wilmington, Del.},
   date={1979},
   pages={xiv+668},
   isbn={0-914098-83-7},
   review={\MR{532830 (82g:53003a)}},
}
\bib{tar}{article}{
   author={Tartar, Luc},
   title={The compensated compactness method applied to systems of
   conservation laws},
   conference={
      title={Systems of nonlinear partial differential equations (Oxford,
      1982)},
   },
   book={
      series={NATO Adv. Sci. Inst. Ser. C Math. Phys. Sci.},
      volume={111},
      publisher={Reidel},
      place={Dordrecht},
   },
   date={1983},
   pages={263--285},
   review={\MR{725524 (85e:35079)}},
}
\bib{te}{article}{
   author={Temple, Blake},
   title={Global solution of the Cauchy problem for a class of $2\times 2$\
   nonstrictly hyperbolic conservation laws},
   journal={Adv. in Appl. Math.},
   volume={3},
   date={1982},
   number={3},
   pages={335--375},
   issn={0196-8858},
   review={\MR{673246 (84f:35091)}},
}
\bib{ts1}{article}{
   author={Tsar{\"e}v, S. P.},
   title={Poisson brackets and one-dimensional Hamiltonian systems of
   hydrodynamic type},
   language={Russian},
   journal={Dokl. Akad. Nauk SSSR},
   volume={282},
   date={1985},
   number={3},
   pages={534--537},
   issn={0002-3264},
   review={\MR{796577 (87b:58030)}},
}
\bib{ts2}{article}{
   author={Tsar{\"e}v, S. P.},
   title={The geometry of Hamiltonian systems of hydrodynamic type. The
   generalized hodograph method},
   language={Russian},
   journal={Izv. Akad. Nauk SSSR Ser. Mat.},
   volume={54},
   date={1990},
   number={5},
   pages={1048--1068},
   issn={0373-2436},
   translation={
      journal={Math. USSR-Izv.},
      volume={37},
      date={1991},
      number={2},
      pages={397--419},
      issn={0025-5726},
   },
   review={\MR{1086085 (92b:58109)}},
}
\bib{yo}{article}{
   author={Young, Robin},
   title={Exact solutions to degenerate conservation laws},
   journal={SIAM J. Math. Anal.},
   volume={30},
   date={1999},
   number={3},
   pages={537--558 (electronic)},
   issn={0036-1410},
   review={\MR{1677943 (2000b:35165)}},
}

\end{biblist}
\end{bibdiv}
\end{document}